\documentclass[final]{siamltex}
\usepackage{amsfonts,amsmath,mathrsfs,bm,mathtools}
\usepackage{amssymb}
\usepackage[notref,notcite]{showkeys}
\usepackage{hyperref}
\usepackage{graphicx}
\usepackage{wrapfig}
\usepackage{subfig}
\usepackage{float}
\usepackage{bbm}
\usepackage{algorithmic, algorithm}
\usepackage{gensymb,color,soul}
\usepackage{tikz}
\usepackage{tikz-3dplot}

\newcommand{\N}{\mathbb{N}}
\newcommand{\R}{\mathbb{R}}
\newcommand{\C}{\mathbb{C}}

\newtheorem{remark}[theorem]{Remark}

\mathtoolsset{showonlyrefs}
\title{Projection-based preprocessing for electrical impedance tomography to reduce the effect of electrode contacts}

\author{
A.~J\"a\"askel\"ainen\footnotemark[2]
\and J.~Toivanen\footnotemark[3]
\and A.~Hänninen\footnotemark[3]
\and V.~Kolehmainen\footnotemark[3]
\and N.~Hyv\"onen\footnotemark[2]
}

\begin{document}
\maketitle

\renewcommand{\thefootnote}{\fnsymbol{footnote}}
\footnotetext[2]{Aalto University, Department of Mathematics and Systems Analysis, P.O.~Box 11100, FI-00076 Aalto, Finland (altti.jaaskelainen@aalto.fi, nuutti.hyvonen@aalto.fi). This work was supported by the Jane and Aatos Erkko Foundation and the Academy of Finland (decision 348503, 353081, 359434, 359181).}
\footnotetext[3]{
University of Eastern Finland, Department of Technical Physics, Kuopio Campus, P.O. Box 1627, FI-70211 Kuopio, Finland (jussi.toivanen@uef.fi, asko.hanninen@uef.fi, ville.kolehmainen@uef.fi). This work was supported by the Jane and Aatos Erkko Foundation and the Academy of Finland (decision 353084, 359433, 358944).}

\begin{abstract}
This work introduces a method for preprocessing measurements of electrical impedance tomography to considerably reduce the effect uncertainties in the electrode contacts have on the reconstruction quality, without a need to explicitly estimate the contacts. The idea is to compute the Jacobian matrix of the forward map with respect to the contact strengths and project the electrode measurements and the forward map onto the orthogonal complement of the range of this Jacobian. Using the smoothened complete electrode model as the forward model, it is demonstrated that inverting the resulting projected equation with respect to only the internal conductivity of the examined body results in good quality reconstructions both when resorting to a single step linearization with a smoothness prior and when combining lagged diffusivity iteration with total variation regularization. The quality of the reconstructions is further improved if the range of the employed projection is also orthogonal to that of the Jacobian with respect to the electrode positions. These results hold even if the projections are formed at internal and contact conductivities that significantly differ from the true ones; it is numerically demonstrated that the orthogonal complement of the range of the contact Jacobian is almost independent of the conductivity parameters at which it is evaluated. In particular, our observations introduce a numerical technique for inferring whether a change in the electrode measurements is caused by a change in the internal conductivity or alterations in the electrode contacts, which has potential applications, e.g., in bedside monitoring of stroke patients. The ideas are tested both on simulated data and on real-world water tank measurements with adjustable contact resistances.
\end{abstract}

\renewcommand{\thefootnote}{\arabic{footnote}}



\pagestyle{myheadings}
\thispagestyle{plain}
\markboth{A.~J\"A\"ASKEL\"AINEN ET AL. 
}{PROJECTIONS TO REDUCE THE EFFECT OF CONTACTS IN EIT}

\section{Introduction}
\label{sec:introduction}

Electrical impedance tomography (EIT) is a method for sensing the internal conductivity of an object by driving electrical currents through it using electrodes placed on the object's surface and measuring resulting voltages on the electrodes. This measurement process defines a nonlinear inverse elliptic boundary value problem for determining the internal conductivity distribution of the object given the measured voltages and the applied currents \cite{Borcea02,Cheney99,Uhlmann09}.

Forming an accurate conductivity reconstruction with EIT can often be difficult due to model uncertainties that may clearly exceed the effect that conductivity changes in the object's interior have on the measurement data~\cite{Barber88,Breckon88,Kolehmainen97}. One of the significant and common sources of error arise from inaccurately known contacts between the electrodes and the surface of the object. The strengths of the contacts can be estimated as auxiliary unknowns (in addition to the conductivity) in the nonlinear reconstruction problem with the cost of a more complicated inverse problem and slightly increased computational cost \cite{Heikkinen02,Vilhunen02}.  Many applications of EIT, such as monitoring of stroke \cite{Toivanen21, Toivanen24} or lung function \cite{Adler09}, are based on dynamic imaging, aiming at the reconstruction of temporal changes in the conductivity using measurements at different time instants. In these applications, estimating the contacts may be nontrivial, especially if employing the widely popular linearized difference imaging algorithms; see,~e.g.,~\cite{Adler09,Barber84}. The linearized approach is often used due to its simplicity, but the reconstructions may suffer from significant artifacts if the contacts change between the measurements, which is a potential hindrance in monitoring of stroke and constitutes the main motivation for this work.

We introduce a method for dealing with the uncertainty in the contacts by projecting the measurement data and the forward map onto the orthogonal complement of the range of the Jacobian matrix of the electrode measurements with respect to the contact strengths. The performance of the approach is further improved if the image space of the projection is also constructed to be orthogonal to the range of the Jacobian with respect to the electrode positions. We show experimentally that while the contact Jacobian depends on the initial estimates for the contacts and the internal conductivity of the imaged object, the subspace onto which one projects is almost independent of them, implying that the approach has high tolerance for uncertainty in the initial estimates of the contact parameters. Although not tested in this work, this observation also motivates forming the projection matrix only once even if resorting to some nonlinear reconstruction method.

The EIT forward problem of solving the electrode potentials given the internal conductivity and applied currents is modelled in this work by the smoothened complete electrode model (CEM) \cite{Hyvonen17b}. The CEM models the surface electrodes as subsets of the examined object's boundary, and it has been found to predict real-world experiments with high accuracy~\cite{Cheng89}. In the smoothened CEM, the contacts are characterized by the contact conductivities over the electrodes, with the option to use smooth contact profiles to achieve better regularity properties for the forward solution. We use the same predefined contact conductivity profile on each electrode and parametrize the strengths of the contacts via their peak values at the centers of the electrodes, which means that the aforementioned contact Jacobian is formed by computing Fr\'echet derivatives with respect to these peak values.

We demonstrate the effects of the projections by computing three-dimensional reconstructions from experimental data measured on a water tank with significant errors in the initial estimates for the contact strengths. Two test cases are considered: the electrode contacts are altered either by partially covering some electrodes by duct tape or by using electrode leads that are equipped with tunable resistors for changing their resistances. Two Bayesian reconstruction algorithms are applied to the forward model linearized with respect to the internal conductivity only, with an accurate initial guess for the conductivity of the salt water filling the tank but highly inaccurate estimates for the contacts on some electrodes. The total variation (TV) prior (cf.~\cite{Rudin92}) combined with the lagged diffusivity iteration \cite{Vogel96} is used to compute reconstructions with clear edges between regions of near-constant conductivity, while the other algorithm resorts to a simple smoothness prior. Although in difference imaging applications the reference data, with respect to which the linearization is performed, are typically {\em measured} at a different time instant (to cancel out modeling errors), we {\em simulate} the reference measurements to demonstrate that our approach also has potential for absolute imaging. For both test cases and algorithms, employing the contact conductivity projection significantly reduces the reconstruction artifacts caused by mismodeling of the contacts, and additionally incorporating the projection with respect to the electrode positions in the algorithms entirely removes the artifacts.

In addition to testing the projection approach in the linearized reconstruction problem of EIT, we also numerically study how the considered projections affect changes in the electrode measurements caused by alterations in the internal conductivity or the contacts. These studies are based on water tank measurements and the three-layer head model introduced in \cite{Candiani19,Candiani22}, and they demonstrate that the signal in the measurements due to altered contacts is significantly reduced by the projections, while the effect on the strength of the signal caused by a change in the interior of the imaged body is much milder. This observation has potential for introducing a test for inferring whether a temporal change in EIT measurements is caused by a change in the conductivity of the brain or an alteration in the contacts in bedside monitoring of stroke patients, cf.~\cite{Toivanen21, Toivanen24}.

This article is structured in the following order. In Section \ref{sec:CEM},  we introduce the smoothened CEM and explain how the derivatives of its solution with respect to different parameters can be computed. Section \ref{sec:comp} describes the computational implementation that is used for the numerical experiments. In Section \ref{sec:data}, we define the projections that are the main topic of this article and numerically study their properties and the effect they have on EIT measurements. The (projected) reconstruction algorithms are introduced in Section \ref{sec:algorithms}, and their performance is investigated in Section \ref{sec:experiments}. Finally, the concluding remarks are presented in Section \ref{sec:conclusion}.

\section{Complete electrode model}
\label{sec:CEM}

In this section, the smoothened version of the CEM is introduced. Additionally, we explain how derivatives of the solution with respect to the conductivity, the contact conductivity and the positions of the electrodes can be computed.

The imaged object is modeled as a bounded Lipschitz domain $\Omega \subset \R^3$ with $M \in \N \setminus \{ 1 \}$ contact electrodes placed onto its surface. We represent the electrodes as open connected subsets $E_1, \dots ,E_M$ of the boundary $\partial \Omega$, on which their closures are mutually disjoint and their union is denoted by $E$. EIT measurements are taken by running a net current pattern $I \in \C^M_\diamond$, where
\[ 
\C^M_\diamond = \Big\{V\in\C^M\,\Big|\, \sum_{m=1}^M V_m = 0\Big\},
\]
through the object via the surface electrodes and measuring the resulting constant electrode potentials. We choose the ground potential level so that the vector of electrode potentials satisfies $U \in \C^M_\diamond$.

The CEM is a mathematical model that can accurately predict real-life EIT measurements \cite{Cheng89}. We use a smoothened version of the CEM, which gives better regularity properties for the electric potential in $\Omega$, while having similar modeling accuracy as the traditional version~\cite{Hyvonen17b}. In the smoothened CEM, the measured potentials $U$ are obtained by finding a weak solution to the elliptic boundary value problem
\begin{equation}
\label{eq:cemeqs}
\begin{array}{ll}
\displaystyle{\nabla \cdot(\sigma\nabla u) = 0 \qquad}  &{\rm in}\;\; \Omega, \\[6pt]
{\displaystyle {\nu\cdot\sigma\nabla u} = \zeta (U - u) } \qquad &{\rm on}\;\; \partial \Omega, \\[2pt]
{\displaystyle \int_{E_m}\nu\cdot\sigma\nabla u\,{\rm d}S} = I_m, \qquad & m=1,\ldots,M. \\[4pt]
\end{array}
\end{equation}
Here $\nu$ denotes the exterior unit normal on $\partial \Omega$. We assume the conductivity $\sigma$ to be isotropic and lie in
\begin{equation}
\label{eq:sigma}
L^\infty_+(\Omega) := \{ \varsigma \in L^\infty(\Omega) \ | \ {\rm ess} \inf {\rm Re} (\varsigma) > 0 \},
\end{equation}
and the contact conductivity $\zeta$ to be a function in
\begin{equation}
  \label{eq:zeta}
\mathcal{Z} := \big\{ \xi \in L^\infty(E) \  \big| \   {\rm Re}\, \xi \geq 0 \   {\rm and} \  {\rm ess} \sup \big( {\rm Re} ( \xi |_{E_m} ) \big) > 0 \  {\rm for} \ {\rm all} \ m= 1, \dots, M \big\}
\end{equation}
that is extended to be a subset of $L^\infty(\partial \Omega)$ using zero-continuation. The functional dependence of the measured potentials $U$ on the parameters is expressed as $U(\sigma,\zeta; I) \in \C^M_\diamond$. As more information on $\sigma$ can be gathered by using up to $M-1$ linearly independent current patterns and increasing the number of current patterns even further can be interpreted as a method to reduce measurement noise, we denote the electrode potentials corresponding to several current patterns as
\begin{equation}
\label{eq:forward_map}
\mathcal{U}\big(\sigma,\zeta; I^{(1)}, \dots, I^{(N)}\big) 
= \big[U(\sigma,\zeta; I^{(1)})^{\top}, \dots, U(\sigma,\zeta; I^{(N)})^{\top}\big]^{\top}
\in \C^{MN},
\end{equation}
where $N \in \N$. The dependence of $\mathcal{U}$ on the employed current patterns is often suppressed, while its dependence on some other parameters such as electrode positions is sometimes written out explicitly.

The problem \eqref{eq:cemeqs} can be expressed as a variational formulation of finding $(u,U) \in \mathcal{H}^1 := H^1(\Omega) \oplus \C^M_\diamond$ satisfying~\cite{Hyvonen17b}
\begin{equation}
\label{eq:weak}
B_{\sigma,\zeta}\big((u,U),(v,V)\big)  \,=  \, I\cdot V \qquad {\rm for} \ {\rm all} \ (v,V) \in  \mathcal{H}^1,
\end{equation}
where $B_{\sigma,\zeta}: \mathcal{H}^1 \times \mathcal{H}^1 \to \C$ is the bounded coercive bilinear form
\begin{equation}
\label{eq:sesqui}
B_{\sigma,\zeta}\big((w,W),(v,V)\big) = \int_\Omega \sigma\nabla w\cdot \nabla v \,{\rm d}x + \int_{\partial \Omega} \zeta (W-w)(V-v)\,{\rm d}S
\end{equation}
and $\cdot$ represents the real dot product. The choice of using a bilinear form and a linear functional in the variational formulation \eqref{eq:weak}, instead of sequilinear and antilinear ones, is intentional since it simplifies certain formulas in the following subsection. Under the presented assumptions, \eqref{eq:weak} has a unique solution in $\mathcal{H}^1$ \cite{Hyvonen17b,Somersalo92}.

\subsection{Fr\'echet derivatives}
In this section, we explain how  Fr\'echet derivatives of the electrode potential component of the solution to the CEM problem with respect to $\sigma$, $\zeta$ and the electrode positions can be computed. The first of these derivatives is needed for forming reconstructions, while the latter two are required for building the respective Jacobians and the associated projections onto the orthogonal complements of their ranges. More detailed information on the Fr\'echet derivatives with respect to a spatially varying, possibly vanishing $\zeta$ and the electrode positions can be found in \cite{Darde12,Darde21,Hyvonen14,Hyvonen17b}; see also \cite{Hyvonen24} that contains a more extensive summary on the needed derivatives in a closely related setting.

The forward solution $U(\sigma,\zeta; I)$ is known to be Fr\'echet differentiable with respect to both $\sigma$ and $\zeta$ \cite{Darde21}. The derivative with respect to  $\sigma$ in the direction of a perturbation $\eta\in L^\infty(\Omega)$, denoted by $D_\sigma U(\sigma; \eta) \in \C_\diamond^M$, can be efficiently built by utilizing the sampling formula
\begin{equation}
\label{eq:sderiv}
D_\sigma U(\sigma; \eta) \cdot \tilde{I} = - \int_{\Omega} \eta \nabla u \cdot \nabla \tilde{u} \, {\rm d} x,
\end{equation}
which contains an auxiliary current pattern $\tilde{I}$. The solution to $\eqref{eq:weak}$ when using $\tilde{I}$ as the current pattern is denoted by $(\tilde{u}, \tilde{U})$. Thus, by evaluating the right-hand side of \eqref{eq:sderiv} for $M-1$ different auxiliary current patterns that form a basis of $\C_\diamond^M$, one can determine the derivative $D_\sigma U(\sigma; \eta)$. Similarly, the Fr\'echet derivative of $U(\sigma,\zeta; I)$ with respect to the contact conductance $\zeta$ can be computed via
\begin{equation}
\label{eq:zderiv}
D_\zeta U(\zeta; \omega) \cdot \tilde{I} = - \int_{\partial \Omega} \omega (U-u) (\tilde{U} - \tilde{u}) \, {\rm d} S,
\end{equation}
where $\omega \in L^\infty(E)$ is the perturbation of $\zeta$.

In order to introduce the derivative with respect to the electrode positions, we assume, for simplicity, that $\partial \Omega$ is of class $\mathcal{C}^\infty$ in some neighborhood of the electrodes and that the electrode boundaries $\partial E$ are also smooth. Moreover, the conductivities $\sigma$ and $\zeta$ are assumed to satisfy $\sigma \in L^\infty_+(\Omega) \cap \mathcal{C}^{0,1}(\overline{\Omega})$ and $\zeta \in \mathcal{Z} \cap H^1(\partial \Omega)$.
Let $a \in \mathcal{C}^1(E, \R^3)$ be a perturbation of the electrodes such that the perturbed versions are defined as
\begin{equation}
\label{eq:perturbed}
E_m^a = \big\{ P_x \big(x +a(x) \big) \, \big| \, x \in E_m \big\} \subset \partial \Omega, \qquad m=1,\dots, M,
\end{equation}
where $P_x: \R^3 \supset B_\rho(x) \to \partial \Omega$ is the projection from a small enough ball of radius $\rho$ around $x$ onto $\partial \Omega$ in the direction of the unit normal $\nu(x)$. The contact conductivity $\zeta^a$ for the perturbed electrodes is defined by
\begin{equation}
\label{eq:perturbedz}
\zeta^a\big(P_x (x + a(x) ) \big) = \zeta(x), \qquad x \in E.
\end{equation}
With this definition, the potentials $U$ are Fr\'echet differentiable with respect to $a$ at the origin, and the formula
\begin{equation}
\label{eq:aderiv}
D_a U(0; h) \cdot \tilde{I} \,=  \,  \int_{\partial \Omega} h_\tau \cdot {\rm Grad}(\zeta) \, (U - u) (\tilde{U}  - \tilde{u}) \, {\rm d} S 
\end{equation}
can be used to compute the associated derivatives, cf.~\cite{Darde12,Hyvonen17b,Hyvonen14}. In \eqref{eq:aderiv}, $h_\tau$ is the tangential component for the direction of movement and ${\rm Grad}(\zeta)$ is the surface gradient of the contact conductivity.

\section{Computational model}
\label{sec:comp}

The computation of the required derivatives and reconstructions is performed by discretizing the domain $\Omega$ using a finite element method (FEM). The finite element (FE) discretization provides a tetrahedral mesh and a set of basis functions $\varphi_j \in H^1(\Omega)$, $j=1, \dots, n$, which can be used to approximate any conductivity distribution as
\begin{equation}
\label{eq:discr_sigma}
\sigma = \sum_{j=1}^n \sigma_j \varphi_j \, .
\end{equation}
The contact conductivity $\zeta$ is represented analogously on the boundary of $\Omega$. In our implementation, we employ piecewise linear ``hat'' basis functions $\{ \varphi_j \}_{j=1}^n$, where $\varphi_j$ attains the value 1 at vertex $j$ of the FE mesh and zero at all other vertices. The solution of \eqref{eq:weak} is approximated in the same FE basis in which the conductivity is presented; we abuse the notation by denoting with $(u,U) \in {\rm span} \{ \varphi_j \}_{j=1}^n  \oplus \C^M_\diamond$ this approximate FE solution in the following. Moreover, any function presented in the FE basis $\{ \varphi_j \}_{j=1}^n$ is identified with the corresponding vector, that is, we may, e.g., treat $u$ either as an element of $H^1(\Omega)$ or as a vector in $\C^n$ depending on which is the more convenient interpretation. For more information on implementing FEM for the CEM, see \cite{Vauhkonen97,Vauhkonen99}.  

The derivatives of the second component of the FE solution $(u,U)$ to \eqref{eq:weak} can be approximated via evaluating the right-hand sides of equations \eqref{eq:sderiv}, \eqref{eq:zderiv} and \eqref{eq:aderiv} for suitable perturbations $\eta$, $\omega$ and $h$. The derivative with respect to $\sigma_j$ is obtained by simply setting $\eta = \varphi_j$ in \eqref{eq:sderiv}, but what is meant by derivatives with respect to the contact conductivity and electrode positions requires further explanations.

All electrodes used in our studies are computationally modeled as disks in the sense that they are defined by intersections of $\partial \Omega$ with a circular cylinder whose radius is $R > 0$ and central axis is parallel to the normal of the surface at the center of the considered electrode. We assume that the contact conductivity $\zeta$ is smooth and takes a fixed form on each electrode, with only its magnitude varying between electrodes.  Restricted to a single electrode $E_m$, $m = 1,\dots,M$, the contact $\zeta$ is thus of the form
\begin{equation}
\label{eq:zeta_param}
\zeta|_{E_m}(r_m, \psi_m) =  \zeta_m \hat{\zeta} (r_m, \psi_m) , \qquad r_m \in [0, R), \ \ \psi_m \in [0, 2 \pi),
\end{equation}
in the polar coordinates $(r_m, \psi_m)$ induced on the electrode patch by the cylinder defining it. We use the infinitely smooth radially symmetric shape function
\begin{equation}
    \hat{\zeta} (r, \psi) = \exp\left( \tau - \frac{\tau R^2}{R^2 - r^2} \right).
\end{equation}
Here $\tau>0$ is a shape parameter, which is set to $\tau = 0.4$ in our numerical experiments. The derivative with respect to the peak contact conductivity $\zeta_m$ can now be obtained by setting $\omega = \hat{\zeta}$ on the $m$th electrode and $\omega = 0$ on the rest of $\partial \Omega$ in \eqref{eq:zderiv}. In what follows, we identify a contact conductivity $\zeta$ with the corresponding vector of peak values in $\C^M$.

Let us then consider the derivative of $U$ with respect to the electrode positions. In all experiments, we define the location of, say, the $m$th electrode using two angles that parametrize its center point: the polar $\theta_m \in [0,\pi/2]$ and the azimuthal angle $\phi_m \in [0,2 \pi)$ with respect to the center of the bottom face of $\Omega$. Movement in the direction of either of these angles is then defined by considering the tangent vectors $\hat{\theta}_m, \hat{\phi}_m \in \R^3$ of $\partial \Omega$ at the electrode center, obtained by differentiating the parametrization of $\partial \Omega$ with respect to the angles $\theta_m$ and $\phi_m$, and extending them to the entire electrode so that the lengths of the respective tangential components are constant over the electrode in order to approximately conserve its area; see \cite[Section~4.3]{Candiani19} for the details. To calculate the derivative with respect to $\theta_m$ or $\phi_m$, we then use
$$
h^\theta_m(x) = \hat{\theta}_m(x)  \quad \text{or} \quad  h^\phi_m(x) =  \hat{\phi}_m(x) ,
\qquad x \in E_m,
$$
as $h$ in \eqref{eq:aderiv}.

Our (numerical) experiments consider a head model adopted from \cite{Candiani19,Candiani22} and experimental water tank data measured by the KIT5 stroke EIT prototype device~\cite{Toivanen21}. These two setups are briefly described in the following two subsections.

\subsection{Computational head model}
\label{sec:head}
\begin{figure}[t]
\center{
  {\includegraphics[width=6.5cm]{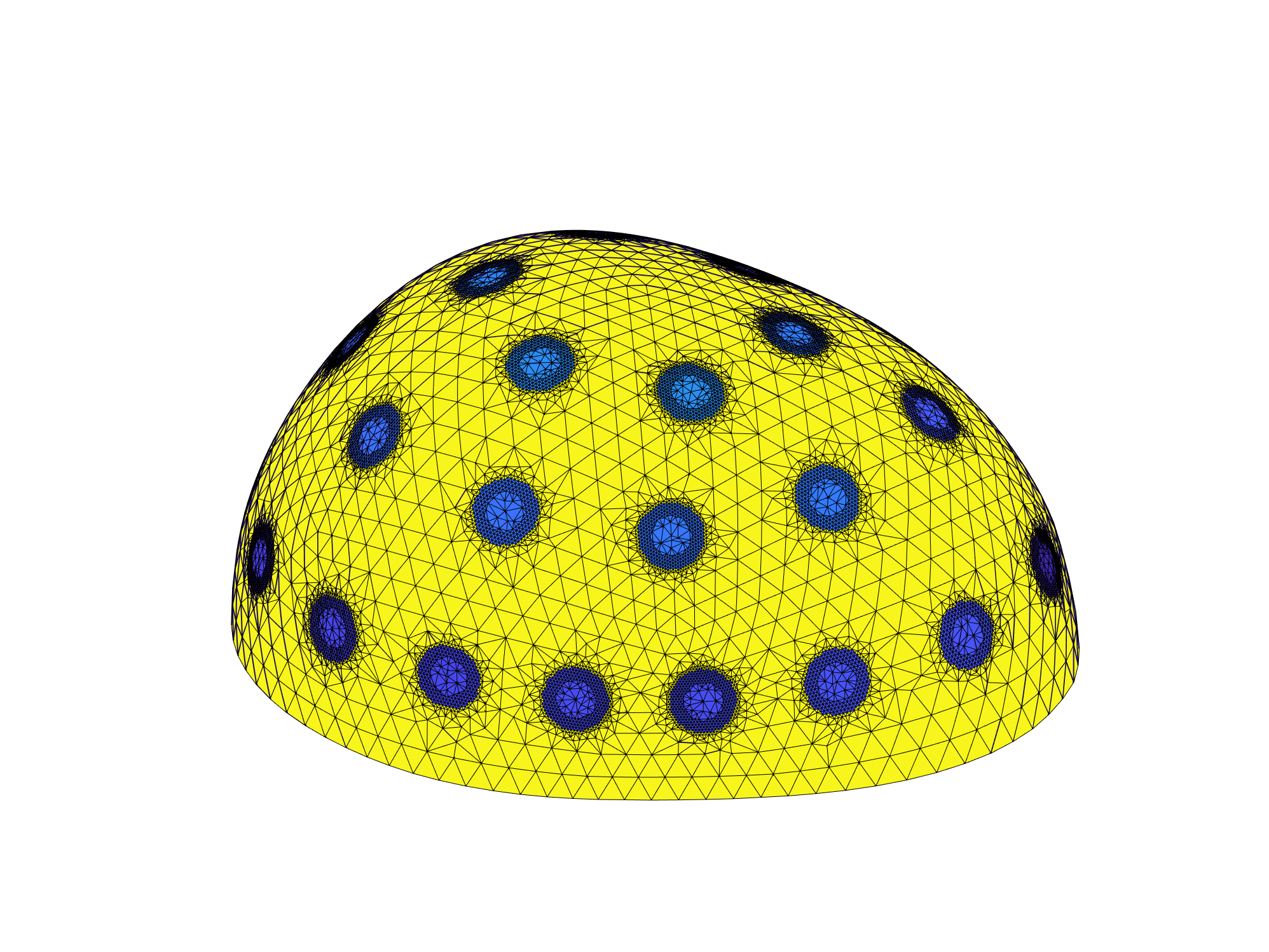}} \quad
  {\includegraphics[width=5.5cm]{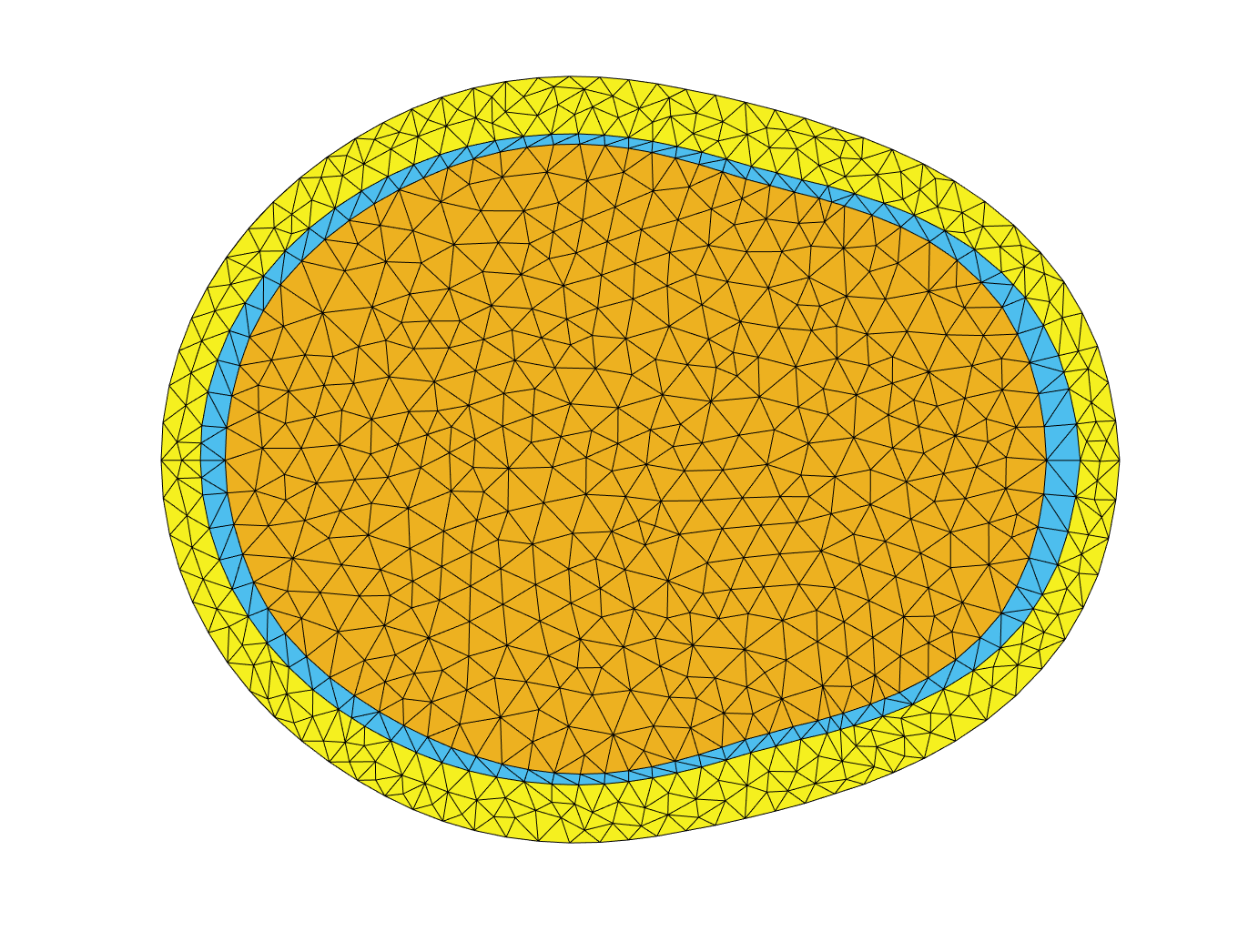}}
  }
\caption{The finite-element head model used for simulating EIT measurements.}
\label{fig:headmodel}
\end{figure}

For the simulated experiments, we use a tetrahedral FE model of a human head that is visualized in Figure \ref{fig:headmodel}. The model consists of three layers: skin, skull and brain. The employed FE mesh has approximately $N = 15\,000$ nodes and $70\,000$ tetrahedra with appropriate refinements at the $M=32$ electrodes of radius $R = 7.5$\,mm distributed approximately uniformly over the head surface. The current patterns are a full basis of $N = M-1$ discrete Fourier currents defined in milliamperes as
$$
I^{(k)} = \left[\cos\left( \frac{2 \pi (m-1) k}{M} \right) \right]_{m=1}^{M}  \quad \text{and} \quad I^{(M/2+l)} =   \left[\sin\left( \frac{2 \pi (m-1) l}{M} \right) \right]_{m=1}^{M},
$$
where $k=1,\dots, M/2$ and $l=1, \dots, M/2-1$. The index $m=1$ corresponds to the frontal electrode on the lowest of the three electrode belts, and the others are belt-wise numbered in the increasing direction of $\phi$ so that the lowest index on each belt is assigned to the frontal electrode. As the base values of the conductivity for the skin, skull and brain we use 0.2\,S/m, 0.06\,S/m and 0.2\,S/m, respectively.

\subsection{Experimental setup}
\label{sec:experimental}
All experimental data were measured using the KIT5 stroke EIT prototype device~\cite{Toivanen21} on a cylindrical measurement tank of radius 11.5\,cm filled with 4.3\,cm of water; see Figure~\ref{fig:taped_electrodes}. The $M=32$ circular electrodes with diameter of 1\,cm on the curved interior surface of the tank were used for current injections of 1\,mA at 12\,kHz frequency, and for measuring the resulting voltages. The first electrode is marked with gray tape in the measurement setup of Figure~\ref{fig:taped_electrodes}, and the others are numbered in the counterclockwise direction. In all experiments, we employ the current patterns ${\rm e}_{2m - 1} - {\rm e}_{2m + 1}$, $m = 1, \dots, 16$, with $\{ {\rm e}_m \}_{m=1}^M$ denoting the Cartesian basis vectors for $\R^M$ and the interpretation ${\rm e}_{33} = {\rm e}_1$. In some tests, the additional current patterns ${\rm e}_{2m - 1} - {\rm e}_{2m + 15}$, $m = 1, \dots, 8$ are also used. As only every second electrode can be used for current input \cite{Toivanen21}, the first set of currents can thus be characterized as {\em adjacent patterns} and the latter as {\em opposite patterns.} The phase information in the electrode voltages is ignored, and their amplitudes are used as measurements for a forward model with real-valued internal and contact conductivities. The cylindrical computational domain, defined by the water layer inside the tank, was discretized into a FE mesh with approximately $N = 25\,000$ nodes and $120\,000$ tetrahedra with appropriate refinements at the electrodes.

\section{Projected data}
\label{sec:data}
In this section, we describe the main idea of this paper, that is, two different projections to remove the effects of uncertain contact conductivities from the potential measurements on the electrodes. In particular, we demonstrate how the range of one of the projections is almost independent of the initially assumed values for the internal and contact conductivities, and we also show examples of this method's effectiveness in recovering a signal caused by a change in the internal conductivity from measurements containing significant modeling error due to unknown or changing contact conductivities. From this point on, we assume that the internal and contact conductivities are real-valued and can thus be identified with vectors in $\R_+^n$ and $\R_+^M$, respectively.

Let us explain how the projections are performed. Following the guidelines in Sections \ref{sec:CEM} and \ref{sec:comp}, we set the stage by computing the derivatives of the electrode potentials $\mathcal{U}(\sigma_0, \zeta_0, \theta_0, \phi_0) \in \R^{MN}$ with respect to the peak contact conductivities and electrode positions to form the respective Jacobian matrices 
$$
J_{\zeta}(\sigma_0, \zeta_0, \theta_0, \phi_0), \qquad J_{\theta}(\sigma_0, \zeta_0, \theta_0, \phi_0), \qquad J_{\phi}(\sigma_0, \zeta_0, \theta_0, \phi_0),
$$
all of which are elements of $\R^{MN \times M}$. Here, $\sigma_0$, $\zeta_0$, $\theta_0$, and $\phi_0$ are, respectively, the assumed values, or initial guesses, for the internal conductivity, the peak contact conductivities and the polar and azimuthal angles of the electrodes. We often suppress the dependence of the Jacobians on these background parameters, especially on $\theta_0$ and $\phi_0$ since the electrode angles are not actively changed in our experiments.

The first, simpler projection is only performed with respect to the effect of the peak contact conductivities and thus only requires the Jacobian matrix $J_{\zeta}$. The orthogonal projection matrix onto the orthogonal complement of the range of $J_{\zeta}$ reads
\begin{equation}
\label{eq:proj}
    P_{\zeta} = \mathrm{I} - J_{\zeta} (J_{\zeta}^{\rm \top} J_{\zeta})^{-1} J_{\zeta}^{\rm \top},
\end{equation}
where $\mathrm{I} \in \R^{MN \times MN}$ is the identity matrix. In all our tests, $J_{\zeta} \in  \R^{MN \times M}$ is not ill-conditioned and has a full rank of $M$, which means that the projection $P_{\zeta}$ can be formed without any numerical complications. 

In order to additionally project with respect to the ranges of $J_{\theta}$ and $J_{\phi}$, a projection matrix, $P_{\zeta, \theta, \phi}$, could again be formed using equation \eqref{eq:proj}, but with $J_{\zeta}$ replaced by a matrix with the combined range of $J_{\zeta}$, $J_{\theta}$ and $J_{\phi}$, i.e., $J_{\zeta, \theta, \phi} = [J_{\zeta} \  J_{\theta} \  J_{\phi}]$. In the experiments that consider projections with respect to electrode positions, however, we only project with respect to the azimuthal angle $\phi$ and therefore use $J_{\zeta, \phi} = [J_{\zeta} \  J_{\phi}]$ in place of $J_{\zeta}$ in \eqref{eq:proj} to form the orthogonal projection $P_{\zeta, \phi} \in \R^{MN \times 2 M}$. As in the case of mere $J_{\zeta}$, the combined Jacobian $J_{\zeta, \phi}$ is well-conditioned and has a full rank of $2M$, meaning that forming $P_{\zeta, \phi}$ is numerically stable. Note that our main aim is {\em not} to use $P_{\zeta, \phi}$ to reduce errors caused by physically misplaced electrodes; it turns out that $P_{\zeta, \phi}$ is better than $P_{\zeta}$ at removing errors caused by mismodeling of the contact strengths in experimental settings, possibly due to some kind of an interplay between contact strengths, positions and shapes.

Recall that the Jacobians $J_{\zeta}(\sigma_0, \zeta_0)$ and $J_{\phi}(\sigma_0, \zeta_0)$ depend on the assumed background parameters $\sigma_0$ and $ \zeta_0$, and thus one would expect the same to be true for the associated orthogonal projections $P_{\zeta}(\sigma_0,\zeta_0)$ and $P_{\zeta, \phi}(\sigma_0,\zeta_0)$ as well. However, our numerical experiments demonstrate that the {\em range} of $P_{\zeta}(\sigma_0,\zeta_0)$ is almost independent of the conductivities $\sigma_0$ and $\zeta_0$. Although the same does not seem to hold for $P_{\zeta, \phi}(\sigma_0,\zeta_0)$ according to numerical experiments not documented here, this observation in any case encourages us to introduce the following approach to solving the reconstruction problem of EIT: only form $P_{\zeta}$ and $P_{\zeta, \phi}$ once at the initial guesses for $\sigma_0$ and $\zeta_0$, multiply the (nonlinear) equation that is to be inverted by one of these matrices, and then proceed to {\em only reconstruct the internal conductivity} of the examined body by considering the projected equation. This approach is demonstrated to produce good quality reconstructions in the numerical experiments of Section~\ref{sec:experiments}, even if the initial guess $\zeta_0$ for the peak contact conductivity values is highly inaccurate.

However, before moving to the actual inverse problem of EIT, the next subsection numerically investigates the angle between the ranges of $P_{\zeta}(\sigma_0, \zeta_0)$ and $P_{\zeta}(\sigma, \zeta)$ with $\sigma$ and $\zeta$ being certain random realizations for the internal and contact conductivities. Subsequently, we study how well the introduced projections are able to remove the effect of inaccuracies in the assumed values for $\zeta_0$, that is, whether it holds that
$$
P_{\zeta} \, \mathcal{U}(\sigma, \zeta) \approx  P_{\zeta} \, \mathcal{U}(\sigma, \zeta_0) \quad \text{or} \quad 
 P_{\zeta, \phi} \, \mathcal{U}(\sigma, \zeta) \approx  P_{\zeta, \phi} \, \mathcal{U}(\sigma, \zeta_0),
$$
where the projections $P_{\zeta}$ and $P_{\zeta, \phi}$ have been computed at the background parameters $\sigma_0$ and $\zeta_0$. If these approximate equalities hold for a wide range of background and perturbed parameters, they can, e.g., be used for testing if a change in electrode measurements in a stroke monitoring application is due to a change in contacts or in the conductivity of the brain.

\subsection{On independence of the range of $P_{\zeta}(\sigma, \zeta)$}
\label{sec:range}

The purpose of this subsection is to numerically verify that the range of the projection $P_{\zeta}(\sigma, \zeta)$ 
is almost independent of the conductivities $\sigma$ and $\zeta$ at which it has been evaluated. We first review what is meant by the principal angles between subspaces following loosely the presentation in~\cite{Bjorck73}. To this end, let $V$ and $W$ be subspaces of $\R^m$ and assume for simplicity that $\dim V = \dim W = q$. Then, the principal angles $\theta_1, \dots, \theta_q$ between $V$ and $W$ are defined recursively via the constrained maximization problems
\begin{align*}
  &\max_{v \in V} \max_{w \in W} v^\top w = v_k^{\top} w_k =: \cos \theta_k, \\
&\| v \|_2 = \| w \|_2 = 1, \\[1mm]
&v^\top v_j = w^\top w_j = 0 \qquad \text{for all} \ \ j=1, \dots, k-1. 
\end{align*}
Note that $\theta_1, \dots, \theta_q$ form an increasing sequence, and hence we call $\theta_{\rm max} := \theta_q$ the largest angle between the subspaces $V$ and $W$.

Consider the head model of Section~\ref{sec:head}. We are interested in the largest angle between the ranges 
\begin{equation}
\label{eq:V_W}
V = \mathcal{R} \big(P_\zeta(\sigma_0, \zeta_0) \big) \quad \text{and} \quad
W = \mathcal{R} \big (P_\zeta(\sigma, \zeta) \big)
\end{equation}
with $m = M (M-1)$ and $q = M(M-2)$. Here, $\sigma_0$ and $\zeta_0$ are the considered background parameters, and $\sigma$ and $\zeta$ are their perturbed versions. Since we expect the angles to be small, it is numerically more stable to compute their sines than cosines \cite{Bjorck73}. It follows from \cite[Theorem~3.4]{Knyazev02} that, in case of \eqref{eq:V_W}, $\sin \theta_q$ is the largest singular value, i.e., the spectral norm, of
$$
P_\zeta(\sigma_0, \zeta_0) \big(\mathrm{I} - P_\zeta(\sigma, \zeta)\big).
$$

The electrode configuration in Figure~\ref{fig:headmodel} defines $\theta_0$ and $\phi_0$, and  $\sigma_0$ corresponds to the conductivity levels for the skin, skull and brain layers defined in Section~\ref{sec:experimental}. The expected peak values of the contact conductivities in $\zeta_0 \in \R^M$ are set to 500\,${\rm S}/{\rm m}^2$, which is approximately in line with water tank experiments with hat-shaped contact conductivities in \cite{Hyvonen17b}. The perturbed conductivities $\sigma$ and $\zeta$ are drawn randomly as follows: The peak contact conductivity on $E_m$ is $[{\rm S}/{\rm m}^2]$
\begin{equation}
\label{eq:rand_zeta}
\zeta_m = 10 + 600 \beta + 380 \upsilon_m, 
\end{equation}
where $\beta$ and $\upsilon_m$ are independent realizations of a uniform random variable on $[0,1]$, with $\beta$ being the same for all electrodes and $\upsilon_m$ redrawn for each $m$. The conductivities of the skin and skull layers in $\sigma$ are assumed to be the same as in $\sigma_0$, but the conductivity of the brain is a realization of a log-normal random variable/field: 
\begin{equation}
\label{eq:rand_sigma}
\sigma_{\rm brain} = \exp{\kappa}.
\end{equation}
In \eqref{eq:rand_sigma}, the exponential function operates componentwise, and $\kappa$ is a realization a random variable whose mean is a constant vector with elements $\log 0.2$ and covariance matrix is defined elementwise as
\begin{equation*}
 \Gamma_{i,j} = \varsigma^2 \exp \left(-\frac{| x_i - x_j |^2}{2\ell^2} \right).
 \end{equation*}
Here, $\ell = 0.02$\,m is the assumed correlation length, $\varsigma = 0.5$ is the pointwise standard deviation, and $x_i$ and $x_j$ are the coordinates of the nodes with indices $i$ and $j$ in the FE mesh of the brain. Note that the mean and standard deviation of any single component in $\sigma_{\rm brain}$ are 0.227\,S/m and 0.121\,S/m, respectively.

In order to test whether $\mathcal{R}(P_\zeta(\sigma, \zeta))$ is independent of the particular realizations of $\sigma$ and $\zeta$, the above procedure of randomly drawing the conductivities was performed 1000 times. Statistics of the largest angle between $\mathcal{R}(P_\zeta(\sigma, \zeta))$ and $\mathcal{R}(P_\zeta(\sigma_0, \zeta_0))$ are presented in Table~\ref{tab:angles}. It is concluded that the angle between the considered ranges is always less than a degree and on average less than half a degree. For comparison, the table also shows the corresponding statistics for the relative discrepancy
$$
{\rm err}_{\rm F} = \frac{\| J_{\zeta}(\sigma,\zeta) - J(\sigma_0,\zeta_0)\|_{\rm F}}{\|J(\sigma_0,\zeta_0)\|_{\rm F}}
$$
between the reference and perturbed Jacobians in the Frobenius norm, demonstrating that the dependence of the elements in $J_{\zeta}(\sigma,\zeta)$ on the pair $(\sigma,\zeta)$ is significant even if that of its range is not.

\begin{table}[t]
    \caption{Statistics for the largest angle between $\mathcal{R}(P_\zeta(\sigma, \zeta))$ and $\mathcal{R}(P_\zeta(\sigma_0, \zeta_0))$ as well as for the relative discrepancy  ${\rm err}_{\rm F}$ between the associated Jacobians over a sample of 1000 random draws of $\sigma$ and $\zeta$ for the head model of Section~\ref{sec:head}.}
    \centering
    \begin{tabular}{ccc ccc}
         $\max \theta_{\rm max} $ & $\mathbb{E} [\theta_{\max}]$ & ${\rm std}[ \theta_{\max}]$ & $\max ({\rm err}_{\rm F}) $ & $\mathbb{E} [{\rm err}_{\rm F}]$ & ${\rm std}[ {\rm err}_{\rm F}]$\\[1mm]
         \hline
         $0.84^{\circ}$ & $0.46^{\circ}$ & $0.12^{\circ}$ & 200.17 & 3.67 & 12.42
    \end{tabular}
    \label{tab:angles}
\end{table}

\subsection{Examples on the action of $P_{\zeta}(\sigma, \zeta)$ and $P_{\zeta,\phi}(\sigma, \zeta)$}
\label{sec:Paction}

Let us assume that the measured electrode potentials $\mathcal{U}(\sigma, \zeta)$ correspond to internal $\sigma$ and contact $\zeta$ conductivities that deviate from the assumed background values $\sigma_0$ and $\zeta_0$. The aim of this section is to demonstrate how the projections $P_\zeta(\sigma_0, \zeta_0)$ and $P_{\zeta, \phi}(\sigma_0, \zeta_0)$ can be used to emphasize the information in $\mathcal{U}(\sigma, \zeta)$ due to the change from $\sigma_0$ to $\sigma$ relative to that due to the change from $\zeta_0$ to $\zeta$. Take note that the projections are always evaluated at the background parameters, i.e., they do not assume any information on $\sigma$ or $\zeta$.

In order to simplify the explanations, we call 
$$
s(\sigma) = \mathcal{U}(\sigma, \zeta_0) - \mathcal{U}(\sigma_0, \zeta_0), \ \ \ s(\zeta) = \mathcal{U}(\sigma_0, \zeta) - \mathcal{U}(\sigma_0, \zeta_0), \ \ \ s(\sigma, \zeta) = \mathcal{U}(\sigma, \zeta) - \mathcal{U}(\sigma_0, \zeta_0),
$$
the $\sigma$-signal, $\zeta$-signal and combined signal, respectively. Denoting either $P = P_\zeta(\sigma_0, \zeta_0)$ or $P = P_{\zeta, \phi}(\sigma_0, \zeta_0)$, it would be optimal if the strength of the projected $\sigma$-signal $P s(\sigma)$ were comparable to $s(\sigma)$, the projected $\zeta$-signal $P s(\zeta)$  almost vanished, and $P s(\sigma, \zeta) \approx P s(\sigma)$. That is, $P$ would not delete from $\mathcal{U}(\sigma, \zeta)$ too much useful information on $\sigma$, but it would project away most of the change caused by the update from $\zeta_0$ to~$\zeta$.

\begin{figure}[t]
\center{
  {\includegraphics[width=5.5cm]{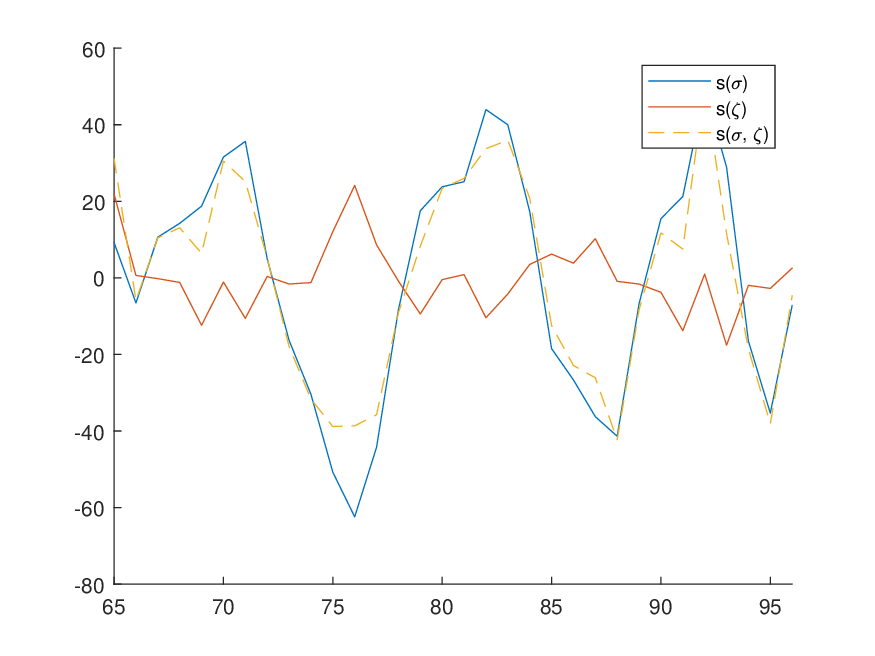}} \quad
  {\includegraphics[width=5.5cm]{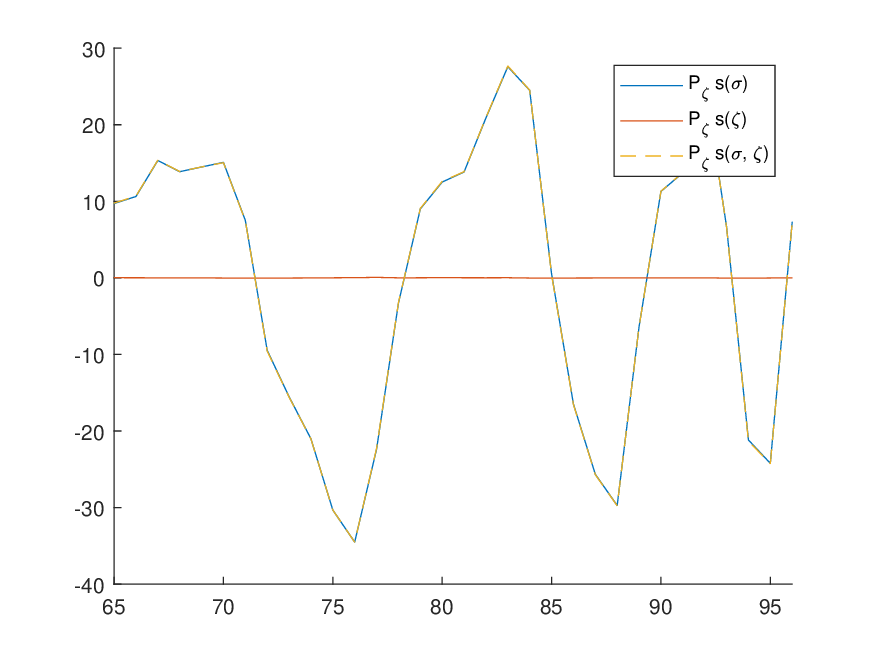}} \\
  {\includegraphics[width=5.5cm]{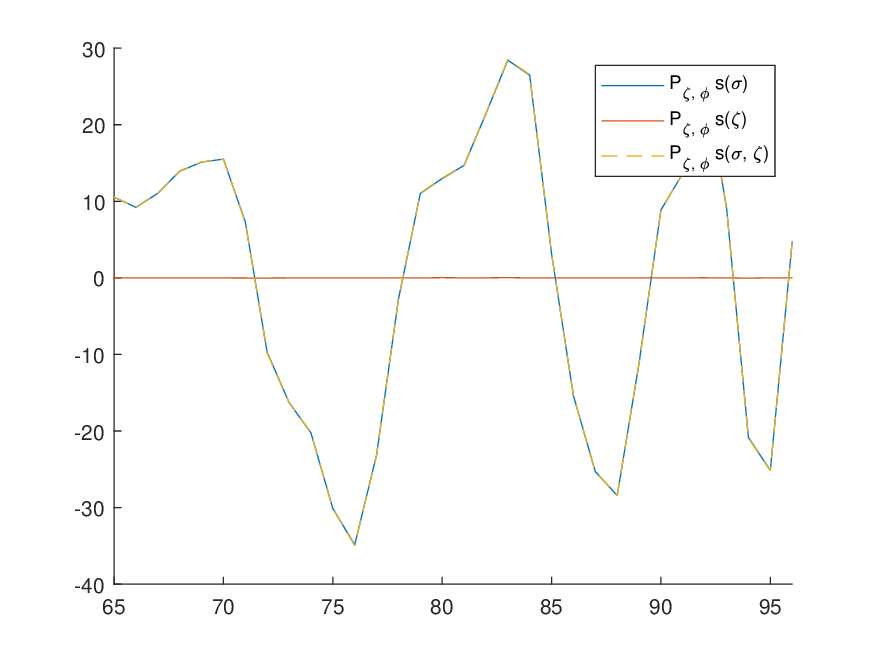}} \quad
  {\includegraphics[width=5.5cm]{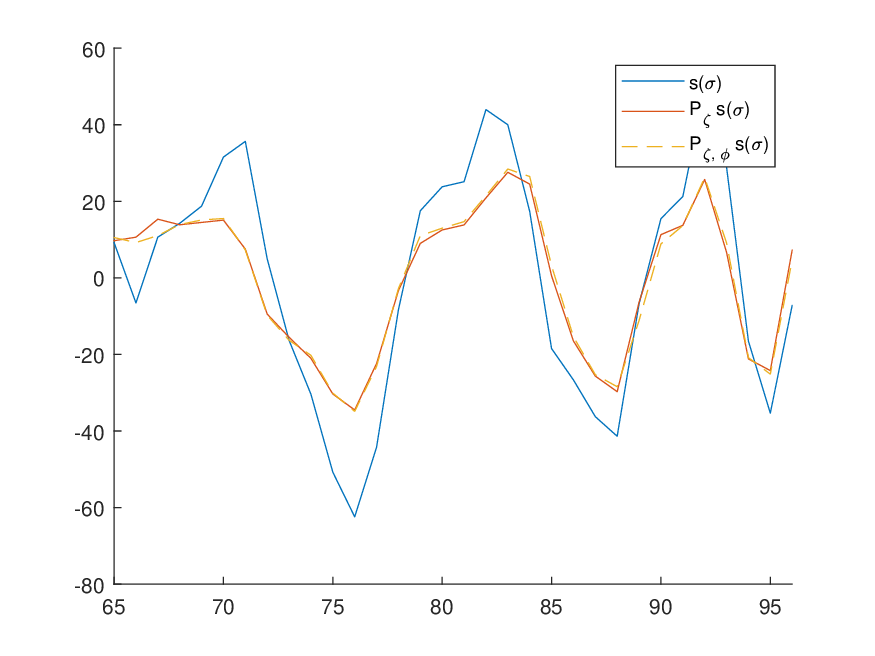}}
  }
\caption{Comparison of $\sigma$-, $\zeta$- and combined signals corresponding to the input current $I^{(3)}$ for a single random draw of the conductivity pair $(\sigma, \zeta)$ in the framework the head model of Section~\ref{sec:head}.
Top left: non-projected signals. Top right: signals projected with $P_\zeta(\sigma_0, \zeta_0)$. Bottom left: signals projected with $P_{\zeta,\phi}(\sigma_0, \zeta_0)$. Bottom right: The original and projected $\sigma$-signals.}
\label{fig:sz_sim}
\end{figure}

We first consider the framework of the head model from Section~\ref{sec:head}. The reference parameters are again chosen so that $\sigma_0$ corresponds to the assumed constant conductivity levels for the skin, skull and brain layers, and all components of $\zeta_0$ equal 500\,${\rm S}/{\rm m}^2$. The perturbed versions $\sigma$ and $\zeta$ are randomly drawn from the same distributions as in Section~\ref{sec:range}.  To first intuitively demonstrate the functioning of $P_\zeta$ and $P_{\zeta,\phi}$, let us consider a single random draw $(\sigma, \zeta)$. Figure~\ref{fig:sz_sim} shows the components of the resulting $\sigma$-, $\zeta$- and combined signals corresponding to the third current pattern $I^{(3)}$, as well as those of their $P_\zeta$ and $P_{\zeta,\phi}$ projected versions.  The projections seem to function almost as desired: although the strength of the $\sigma$-signal is undesirably weakened by the projections, on the positive side the projected $\zeta$-signals almost vanish, and the projected $\sigma$- and combined signals nearly coincide. In particular, the difference between the effects of $P_\zeta$ and $P_{\zeta,\phi}$ is almost negligible in this test based on simulated data.

To confirm the conclusion made based on Figure~\ref{fig:sz_sim}, Table~\ref{tab:signal_norms_sz} shows the mean values for the $2$-norms of $s(\sigma)$, $s(\zeta)$ and $s(\sigma, \zeta)$ as well as those for the corresponding projected signals over 1000 random draws of the perturbed conductivity pair $(\sigma, \zeta)$. Although both projections considerably decrease the norm of the $\sigma$-signal ($P_{\zeta}$ by $50.7$\% and $P_{\zeta,\phi}$ by $52.3$\% on average), the corresponding numbers for the $\zeta$-signal are much higher, namely $99.7$\% and $99.8$\%, respectively. Moreover, the norms of $P_{\zeta} s(\sigma)$ and $P_{\zeta,\phi} s(\sigma)$ are, respectively, on average only $0.2$\% smaller than those of $P_{\zeta} s(\sigma,\zeta)$ and $P_{\zeta,\phi} s(\sigma,\zeta)$, even though the corresponding relative difference for the non-projected signals $s(\sigma)$ and $s(\sigma,\zeta)$ is $11$\%.

\begin{figure}[h!]
\center{
  \includegraphics[width=4.2cm]{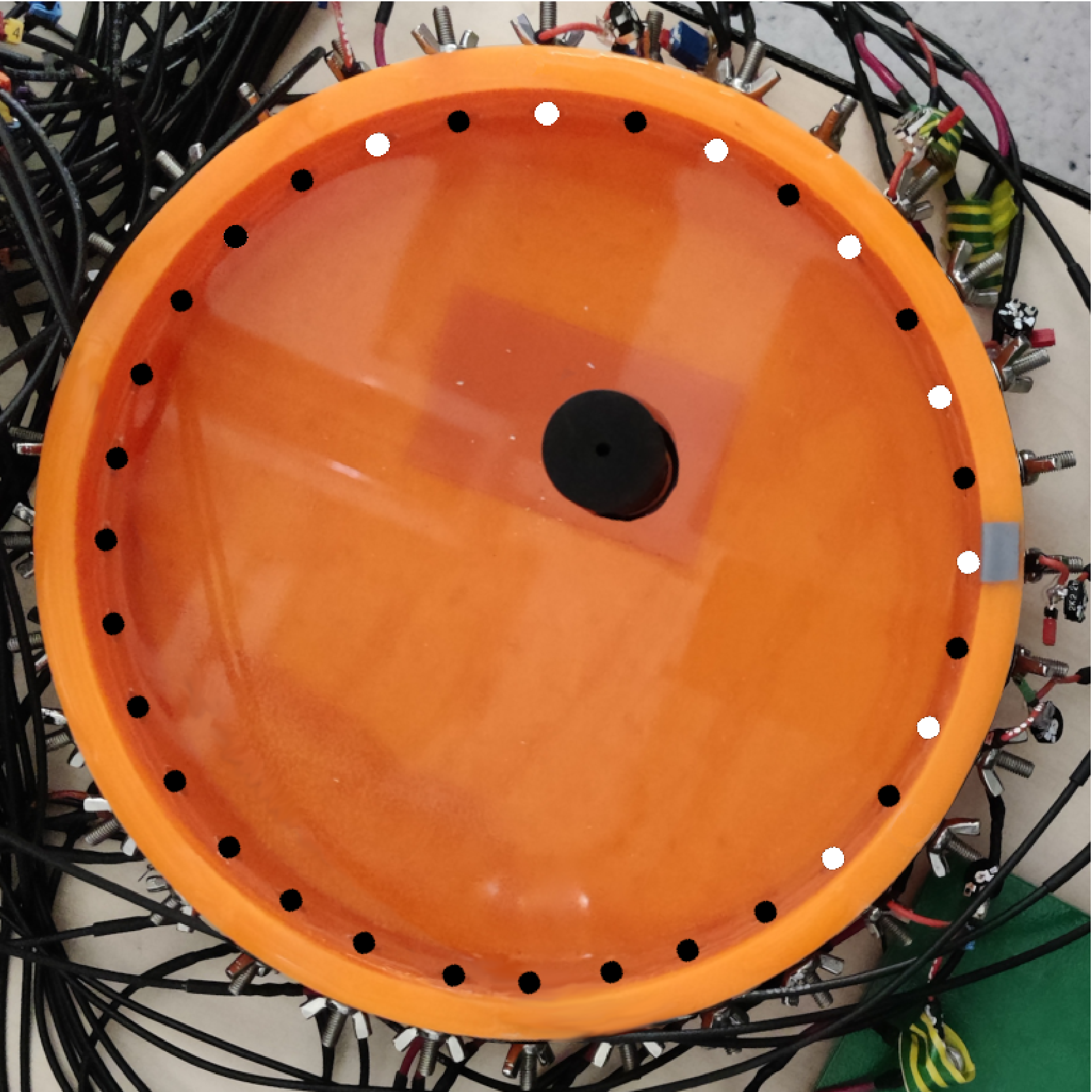} \qquad \qquad
  \includegraphics[width=4.2cm]{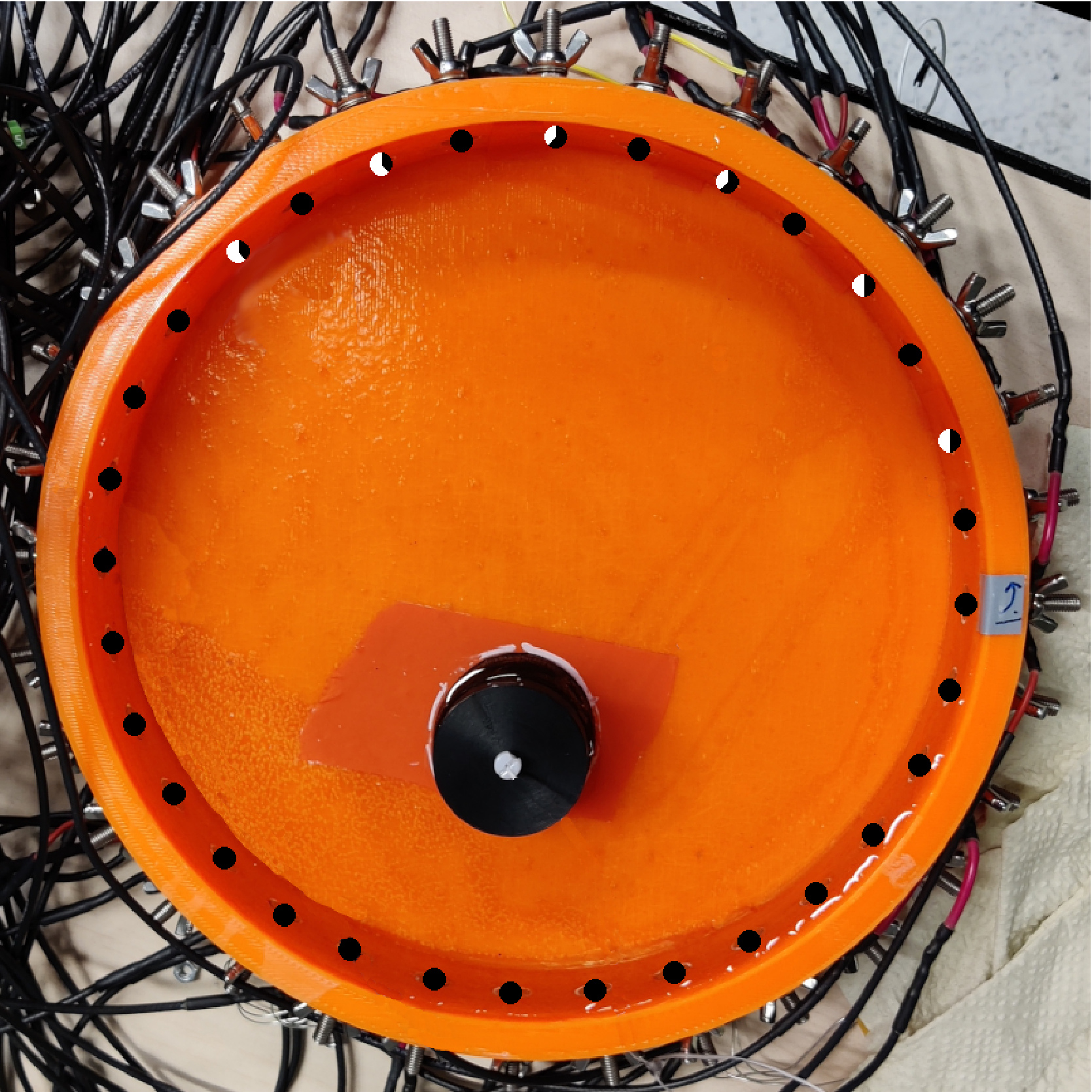}
  }
\caption{Experimental setup described in Section~\ref{sec:experimental}, with the resistor test case on left and the tape test case on right. Left: the electrodes with adjustable resistors in their cables are marked in white. Right: the taped electrodes are partially white, with a half, one third or two thirds of their respective areas covered by tape. See Table~\ref{tab:resistors_and_tapes} for more details. 
}
\label{fig:taped_electrodes}
\end{figure}

These computational results are further confirmed by two tests in the experimental setup described in Section~\ref{sec:experimental}. In both tests, four measurements were taken: (i) with an empty tank, i.e., with the tank filled with salt water only, (ii) with a cylindrical inclusion embedded in the tank and no alterations to the contacts, (iii) with the electrode contacts considerably worsened at some electrodes and no change to the internal conductivity, and (iv) with both an embedded inclusion and altered contacts. The two tests differ in respect to the method used for changing the contacts: in the {\em resistor test case}, changes in contact resistances were mimicked by using adjustable resistors in the electrode cables, whereas in the {\em tape test case}, the areas of some electrodes were changed by partially covering them with tape. Both adjacent and opposite current patterns were used in the resistor test case, whereas the tape test case only considered adjacent patterns (see Section~\ref{sec:experimental}). The resistance values of the adjustable resistors used in the resistor test case and the amount of tape coverage on different electrodes in the tape test case are presented in Table \ref{tab:resistors_and_tapes}. Photos of the corresponding setups, with modified electrodes highlighted, are shown in Figure \ref{fig:taped_electrodes}.

\begin{table}
  \caption{Sample means for the 2-norms of the $\sigma$-, $\zeta$- and combined signals as well as their $P_\zeta$ and $P_{\zeta,\phi}$ projected versions over 1000 random draws of the conductivity pair $(\sigma, \zeta)$ in the framework of the head model of Section~\ref{sec:head}.
    }
    \centering
    \begin{tabular}{c|ccc}
         Projection\textbackslash Signal & $s(\sigma)$ & $s(\zeta)$ & $s(\sigma, \zeta)$ \\
             \hline
         No projection & 656.58 & 368.57 & 737.26 \\
         $P_\zeta(\sigma_0, \zeta_0)$ & 323.71 & 1.23 & 324.47 \\
         $P_{\zeta, \phi}(\sigma_0, \zeta_0)$ & 313.22 & 0.90 & 313.93
    \end{tabular}
    \label{tab:signal_norms_sz}
\end{table}

%

\begin{table}[t]
  \caption{Resistance values for the resistor test case and the tape coverage for the tape test case in the framework of water tank experiments of Section~\ref{sec:experimental}; see Figure~\ref{fig:taped_electrodes}.
        }
    \centering
    \begin{tabular}{c|ccccccccc}
         Electrode & 1 & 3 & 5 & 7 & 9 & 11 & 13 & 29 & 31 \\
            \hline
         R [$\Omega$] & 2200 & 2700 & 2600 & 1000 & 800 & 2000 & - & 1900 & 2400 \\
         Tape [\%] & - & 50 & 50 & 33 & 33 & 66 & 66 & - & - 
    \end{tabular}
    \label{tab:resistors_and_tapes}
\end{table}

In both test cases, the used inclusion was a right circular cylinder with conductivity 4.73\,S/m and height 54\,mm. The radius of the cylinder was 15\,mm in the resistor test case and 20\,mm in the tape test case. The projections $P_\zeta = P_{\zeta}(\sigma_0, \zeta_0)$ and $P_{\zeta,\phi}= P_{\zeta, \phi}(\sigma_0, \zeta_0)$ were evaluated at the measured conductivity of the salt water filling the tank (i.e., $\sigma_0 = 0.0491$\,S/m for the resistor case and $\sigma_0 = 0.0322$\,S/m for the tape case) and for peak contact conductivity values estimated for the empty tank with no changes to the electrode contacts; see, e.g., \cite[Section~4.3]{Hyvonen17b} for details on how the peak values can be estimated. However, the precise values for $\sigma_0$ and $\zeta_0$ do not have a major effect on the presented results. 

\begin{figure}[t]
\center{
  {\includegraphics[width=5.5cm]{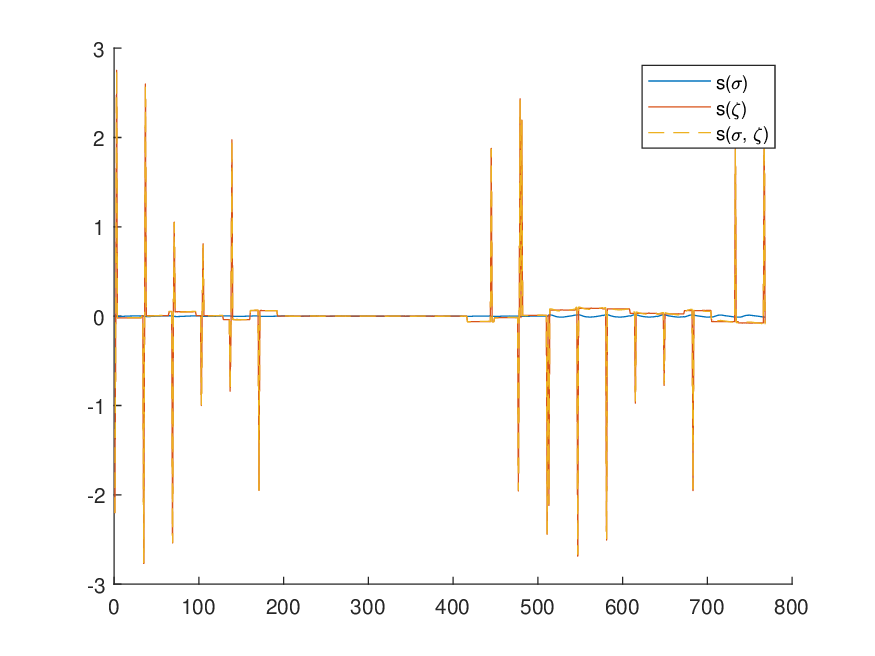}} \quad
  {\includegraphics[width=5.5cm]{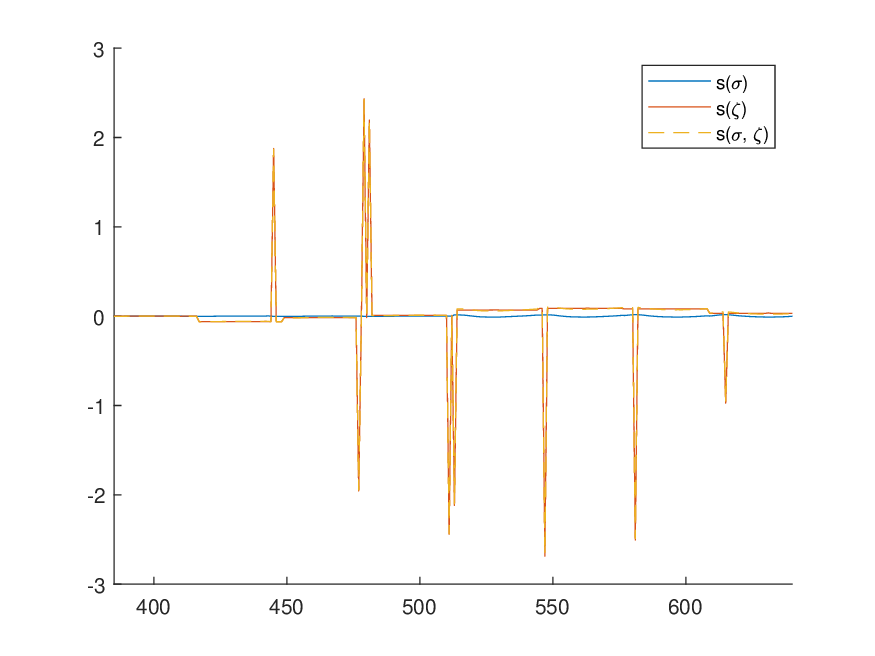}}
  \linebreak
  {\includegraphics[width=5.5cm]{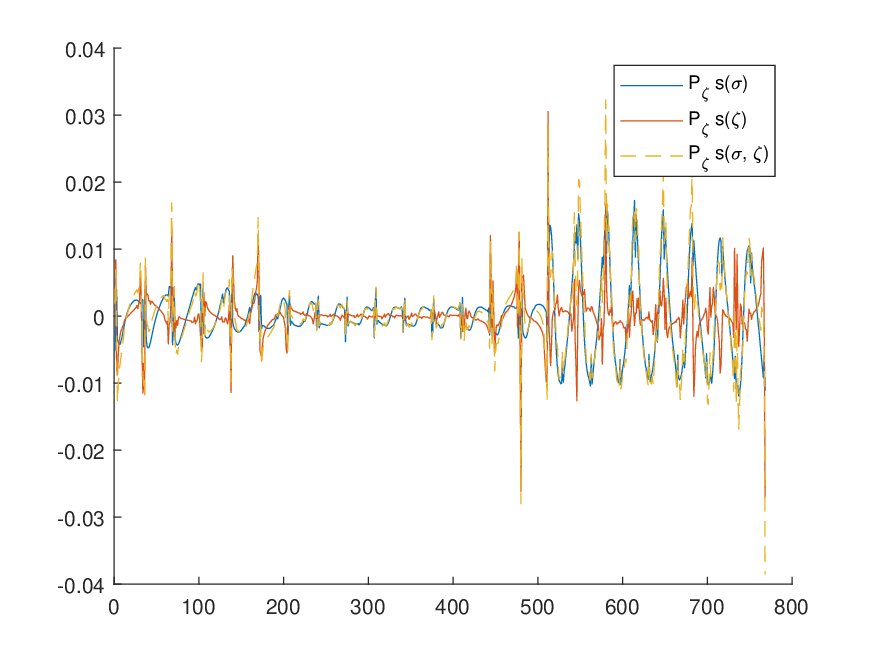}} \quad
  {\includegraphics[width=5.5cm]{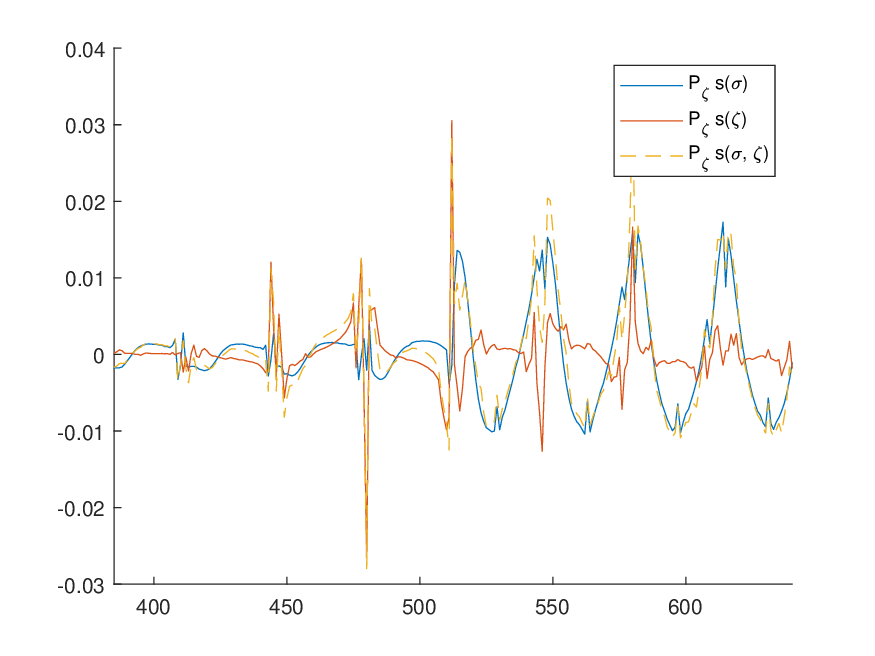}}
  \linebreak
  {\includegraphics[width=5.5cm]{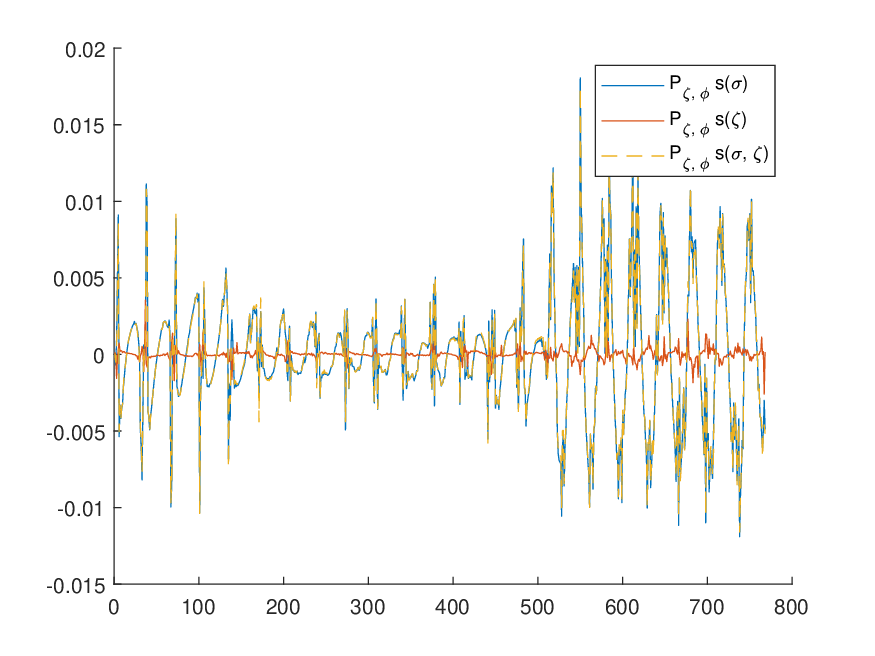}} \quad
  {\includegraphics[width=5.5cm]{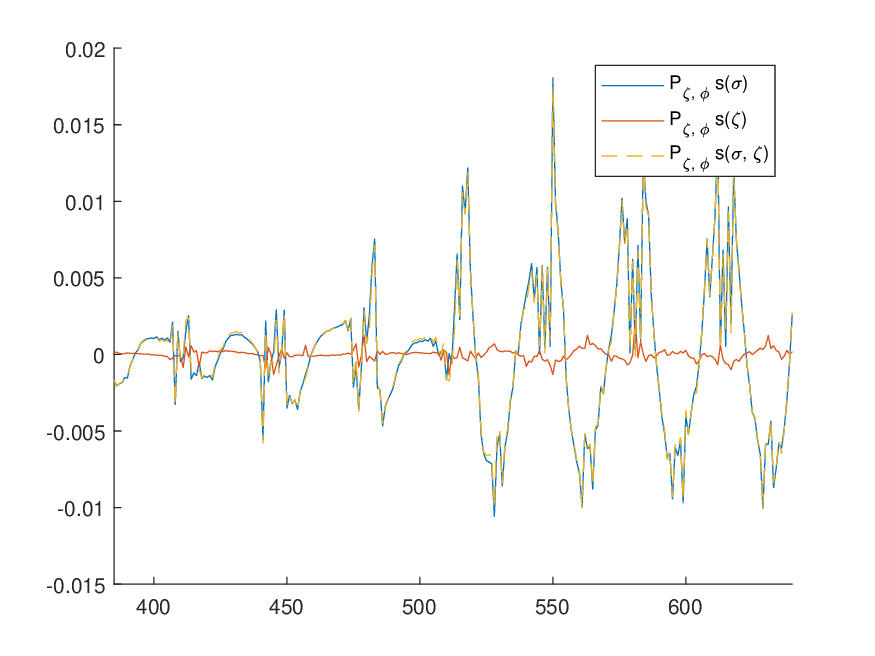}}
  }
\caption{Comparison of experimentally measured $\sigma$-, $\zeta$- and combined signals in the resistor test case; see the left image in Figure~\ref{fig:taped_electrodes}.
 Top: non-projected signals. Middle: signals projected with $P_\zeta(\sigma_0, \zeta_0)$. Bottom: signals projected with $P_{\zeta, \phi}(\sigma_0, \zeta_0)$. Left: full signals corresponding to measurements on all electrodes. Right: zoomed in image on a portion of the signals.}
\label{fig:R_data}
\end{figure}

\begin{figure}[t]
\center{
  {\includegraphics[width=5.5cm]{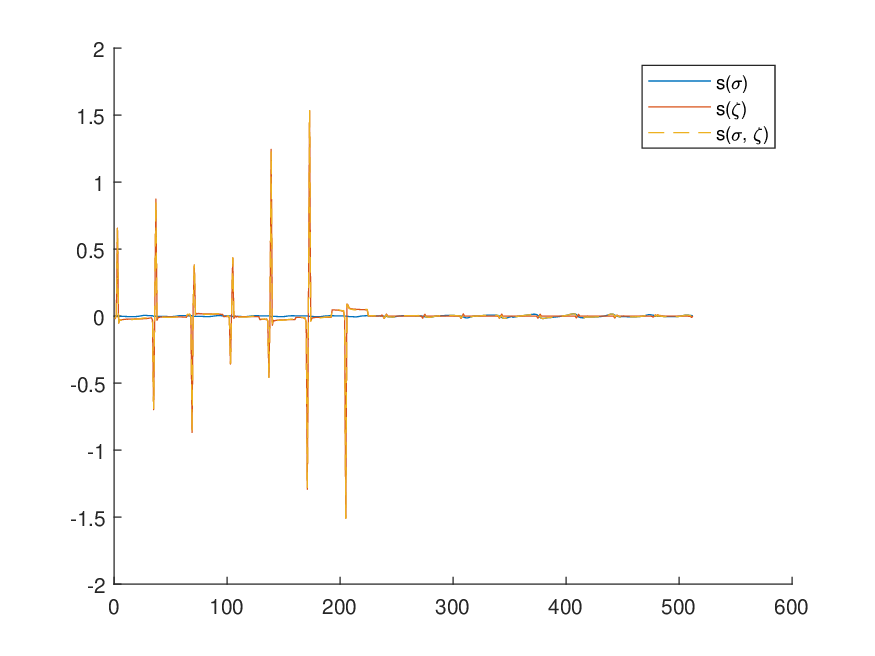}} \quad
  {\includegraphics[width=5.5cm]{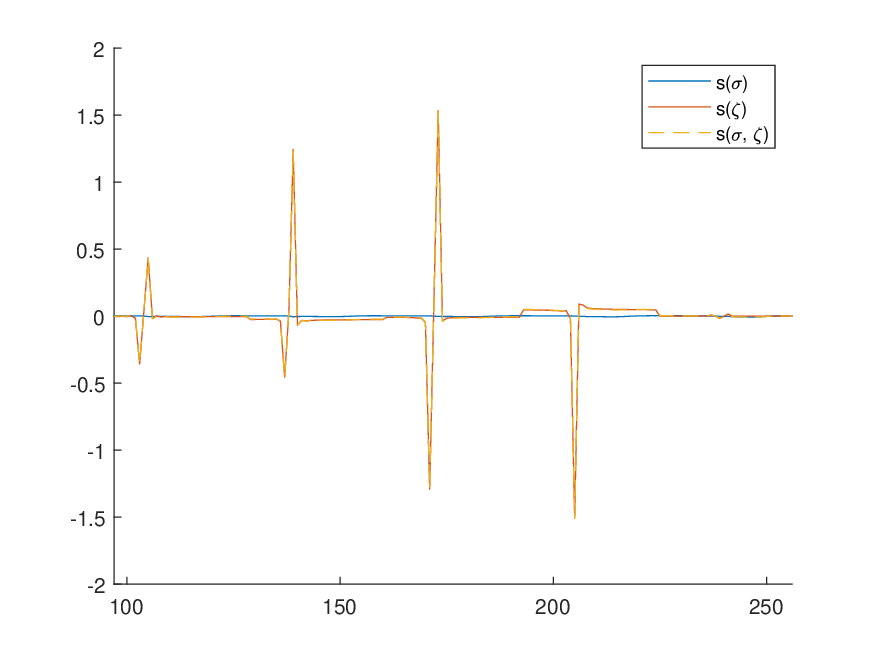}}
  {\includegraphics[width=5.5cm]{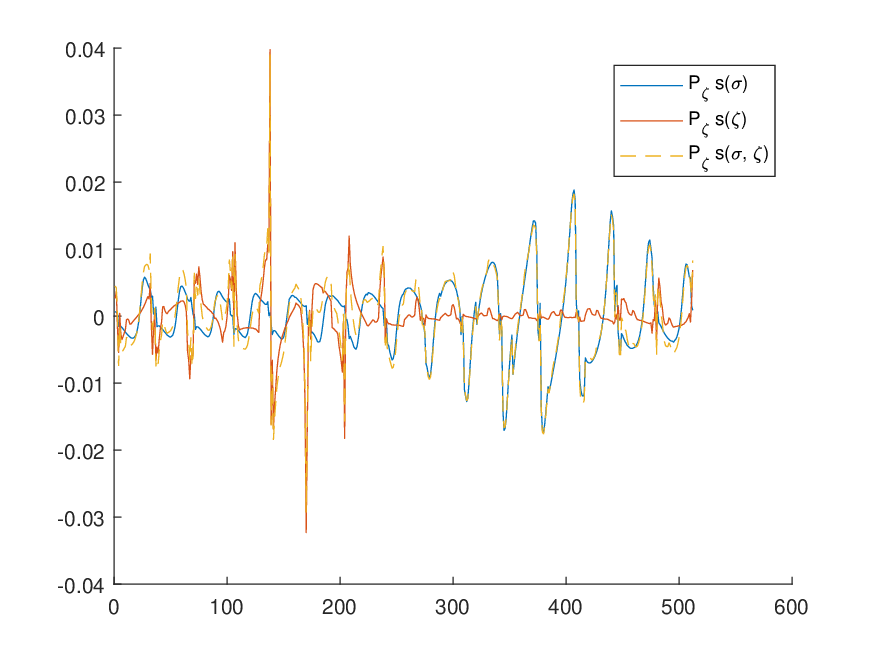}} \quad
  {\includegraphics[width=5.5cm]{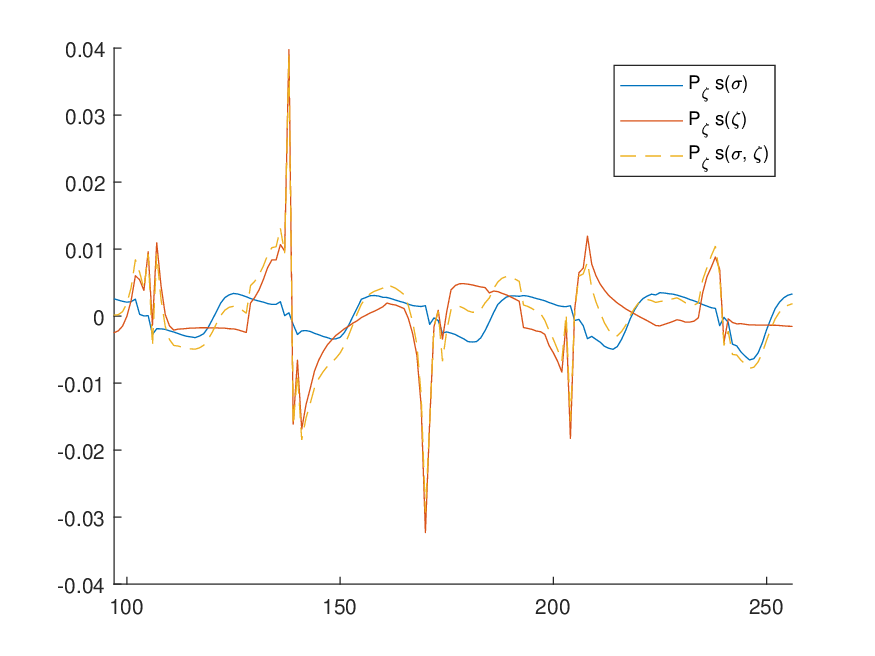}}
  {\includegraphics[width=5.5cm]{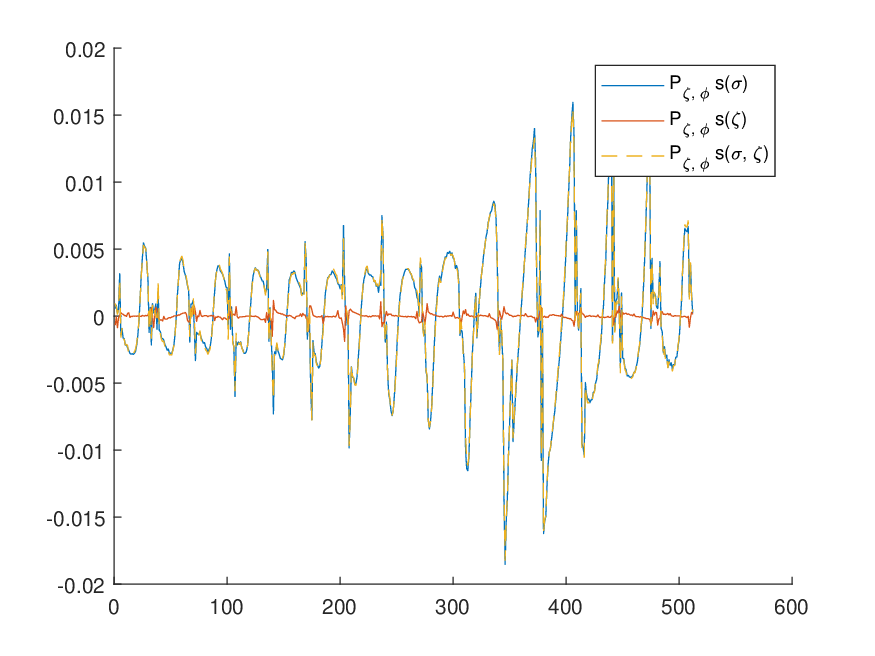}} \quad
  {\includegraphics[width=5.5cm]{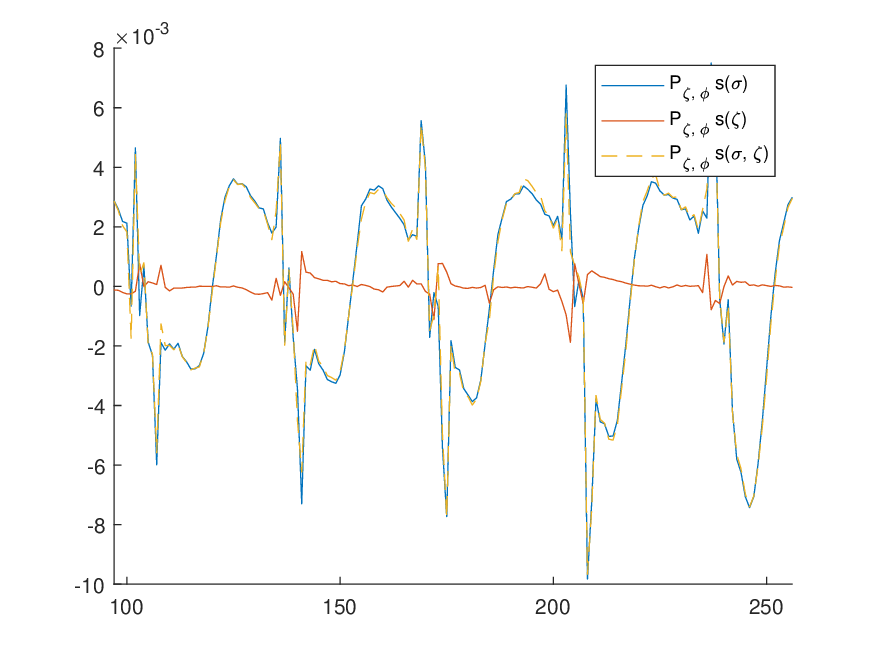}}
  }
\caption{Comparison of experimentally measured $\sigma$-, $\zeta$- and combined signals in the tape test case; see the right image in Figure~\ref{fig:taped_electrodes}.
 Top: non-projected signals. Middle: signals projected with $P_\zeta(\sigma_0, \zeta_0)$. Bottom: signals projected with $P_{\zeta, \phi}(\sigma_0, \zeta_0)$. Left: full signals corresponding to measurements on all electrodes. Right: zoomed in image on a portion of the signals.}
\label{fig:tape_data}
\end{figure}

The results are documented in Figures~\ref{fig:R_data} and \ref{fig:tape_data}. The top rows correspond to the non-projected signals $s(\sigma)$, $s(\zeta)$ and $s(\sigma, \zeta)$, the second row to their $P_{\zeta}$ projected versions and the bottom rows to the $P_{\zeta, \phi}$ projected versions to account also for a possible change in the shape and azimuthal angle of the actual contact. The left column in both figures shows all components of the respective signals, whereas the right column highlights certain components to enable visual evaluation of finer details. For both the resistor test case and the tape test case, the signal caused by the changes in the contacts completely dominates that originating from the embedded inclusion; on the top rows of Figures~\ref{fig:R_data} and \ref{fig:tape_data}, the signals $s(\zeta)$ and $s(\sigma, \zeta)$ almost coincide, while $s(\sigma)$ caused by the inclusion is hardly visible. Projecting with $P_{\zeta}$ is in both test cases able to significantly reduce $s(\zeta)$, but one cannot claim that the projected signals $P_{\zeta}s(\sigma)$ and $P_{\zeta}s(\sigma,\zeta)$ coincide, that is, there is still a clear contribution from the change in the contacts in $P_{\zeta}s(\sigma,\zeta)$. This is not surprising for the tape test case since covering part of an electrode definitely moves the contact area, but based on the second row in Figure~\ref{fig:R_data}, $P_{\zeta}$ cannot completely project away the contribution of $\zeta$ to $s(\sigma,\zeta)$ for the resistor case either. However, as depicted on the last rows of Figures~\ref{fig:R_data} and \ref{fig:tape_data}, the $P_{\zeta, \phi}$ projected signals $P_{\zeta,\phi}s(\sigma)$ and $P_{\zeta, \phi} s(\sigma,\zeta)$ almost exactly coincide for both test cases, that is, the effect of $\zeta$ is almost invisible in $P_{\zeta, \phi}s(\sigma,\zeta)$. Moreover, the amplitude of $P_{\zeta, \phi}s(\zeta)$ is low, meaning that $P_{\zeta, \phi}$ also seems to be able to reveal if the observed change in the measurements originates from worsened contacts only. Although not easily visible in Figures~\ref{fig:R_data} and \ref{fig:tape_data}, it should also be noted that the strengths of the original and projected $\sigma$-signals, $s(\sigma)$ and $P_{\zeta, \phi}s(\sigma)$, are almost the same in both cases.

\section{Projected reconstruction algorithms}
\label{sec:algorithms}

This section describes two algorithms that are used to compute reconstructions of the (change in the) conductivity and to compare the results between cases when the projections are used and when they are not. The algorithms only utilize absolute measurements from the experimental setup described in Section~\ref{sec:experimental}, i.e., no reference data measured from the empty tank are used as their input. However, we do assume to know the measured conductivity of the salt water filling the tank. Moreover, when computing reference reconstructions with ``accurate" information on the contacts, we utilize the peak contact conductivity values that are estimated based on empty tank measurements with no modifications to the contacts. These values and the background conductivity level of the salt water are also used for forming the projections $P_\zeta$ and $P_{\zeta, \phi}$. In particular, the reconstruction algorithms do not estimate the contacts, and apart from the reference reconstructions, the utilized information on the peak contact conductivity values at some electrodes is highly inaccurate since the electrode contacts are deliberately worsened as described in Section~\ref{sec:Paction}.

In their basic forms, both algorithms are based on a single linearization of the forward model with respect to the internal conductivity of the imaged body; see~Remark~\ref{remark:multilin} below for generalizations. To this end, set $y = \mathcal{U}(\sigma, \zeta) - \mathcal{U}(\sigma_0, \zeta_0)$ and consider the Jacobian $J_\sigma = J_\sigma(\sigma_0, \zeta_0)$ (cf.~Section \ref{sec:CEM}), with $(\sigma, \zeta)$ denoting the true conductivity pair characterizing the object and the contacts, and $(\sigma_0, \zeta_0)$ corresponding to the best available information on them. In our setting, $\mathcal{U}(\sigma, \zeta)$ is measured and $\mathcal{U}(\sigma_0, \zeta_0)$ computed, with $\sigma_0$ being a good estimate for the true conductivity away from the inhomogeneity and $\zeta_0$ typically being a significantly inaccurate estimate for $\zeta$ on some electrodes.
Both algorithms are based on a Bayesian approach, and thus we model the measurement data $y$ as a realization of a random variable that is defined by the formula
\begin{equation}
\label{eq:lin_mod}
    Y = J_\sigma W + E.
\end{equation}
Here, the random variable $W$ models the difference of the true conductivity from the assumed background value, i.e., it corresponds to a randomization of $\sigma - \sigma_0$. The noise vector $E$ is a zero-mean Gaussian random variable that is assumed to be independent of the change in the conductivity $W$. The projected version of \eqref{eq:lin_mod} reads
\begin{equation}
\label{eq:proj_lin_mod}
    P Y = P J_\sigma W + P E,
\end{equation}
where $P = P_{\zeta}(\sigma_0, \zeta_0)$ or $P = P_{\zeta, \phi}(\sigma_0, \zeta_0)$.

It is assumed that $W$ follows a smoothened TV-type prior (cf.~\cite{Rudin92}) defined by 
\begin{equation}
  \label{eq:TV}
\pi(w) \propto \exp (- \gamma \Psi(w)  ) , \qquad \gamma > 0,
\end{equation}
 where
\begin{equation}
\label{eq:TVexp}
\Psi(w) = \int_{\Omega} \psi \big(|\nabla w | \big) \, {\rm d} x + \frac{\varepsilon}{2}  \|w\|_2^2, \qquad \text{with} \quad \psi(t) = \sqrt{t^2 + T^2} \approx | t |,
\end{equation}
and $T, \epsilon > 0$ are {\em small} parameters ensuring, respectively, that $\psi$ is differentiable and \eqref{eq:TV} defines a proper prior. To understand \eqref{eq:TVexp}, recall that $w$ is interpreted as an element of $H^1(\Omega)$ via the employed FE basis when appropriate (cf.~\eqref{eq:discr_sigma}). Such a prior is known to promote monotonic changes and prevent significant amounts of oscillations in the reconstruction.

Denote by $\Gamma_E$ the symmetric positive definite covariance matrix of $E$. The MAP estimates for \eqref{eq:lin_mod} and \eqref{eq:proj_lin_mod}, i.e., the points of highest posterior probability densities defined by \eqref{eq:lin_mod} and \eqref{eq:proj_lin_mod} together with \eqref{eq:TV}, are characterized by minimizers of
\begin{equation}
\label{eq:Tikhonovk}
 \frac{1}{2} (y - J_\sigma w )^{\rm \top} \Sigma (y - J_\sigma w ) +  \gamma \, \Psi(w) 
\end{equation}
where $\Sigma = \Gamma_E^{-1}$ for  \eqref{eq:lin_mod} and $\Sigma = P \Gamma_E^{-1} P$ for \eqref{eq:proj_lin_mod}.
To deduce the latter, note that the random vector $PE$ follows on the subspace $\mathcal{R}(P)$ a Gaussian distribution, which can be obtained, e.g., by marginalizing that of $E$ over $\mathcal{R}(P)^\perp$, with the covariance matrix $P \Gamma_E P^\top = P \Gamma_E P$. The inverse of this covariance matrix on $\mathcal{R}(P)$ is realized by the pseudoinverse of $P \Gamma_E P$, i.e., $P \Gamma_E^{-1} P$. The formula \eqref{eq:Tikhonovk} in case of \eqref{eq:proj_lin_mod} thus follows by combining the likelihood function induced by \eqref{eq:proj_lin_mod} with the prior \eqref{eq:TV} and noting that for an orthogonal projector $P^2 = P = P^\top$.

Our first, simpler algorithm corresponds to using a certain quadratic approximation for $\Psi(w)$ around $w=0$ in \eqref{eq:Tikhonovk}, which leads to standard Tikhonov regularization with a smoothness prior for the linearized (and possibly projected) equation. The second algorithm finds a minimizer for \eqref{eq:Tikhonovk} via the lagged diffusivity iteration.

\subsection{Regularized one-step linearization}
\label{sec:one-step}
 To tackle minimizing \eqref{eq:Tikhonovk}, we consider the necessary optimality condition
\begin{equation}
\label{eq:necessary}
J_\sigma^{\rm \top} \Sigma J_\sigma  w  
+ \gamma (\nabla_w \Psi)(w) 
 =  J_\sigma^{\rm \top} \Sigma y.
\end{equation}
The gradient term in \eqref{eq:necessary} can be represented as \cite{Arridge13}
\begin{equation}
  \label{eq:kappa_deriv}
(\nabla_w \Psi)(w) = \Theta(w) w,
\end{equation}
where the matrix $\Theta(w)$ is defined elementwise by
$$
\Theta_{i,j}(w) := \int_{\Omega} \frac{\nabla \varphi_i(x) \cdot \nabla \varphi_j(x)}{\sqrt{|\nabla_x w (x)|^2 + T^2}} \, {\rm d} x + \varepsilon \delta_{i,j} , \qquad i,j=1,\dots ,n,
$$
with $\delta_{i,j}$ denoting the Kronecker delta. Due to an interpretation as the system matrix for a FE discretization of an elliptic operator with a homogeneous Neumann boundary condition, with its spectrum shifted by $\epsilon$, the matrix $\Theta(w)$ is invertible, cf.~\cite{Harhanen15}.

The idea of the first algorithm is to replace $\Theta(w)$ by $\Theta(0)$ in \eqref{eq:kappa_deriv} and use the resulting approximation for $(\nabla_w \Psi)(w)$ in  \eqref{eq:necessary}, which leads to a linear equation. This procedure can be interpreted as adopting a certain quadratic approximation for $\Psi$ around the origin \cite{Helin22}. Introducing in addition a symmetric factorization $\Sigma = B^\top B$ and setting 
\begin{equation}
\label{eq:Amatrix}
A = B J_\sigma \quad \text{and} \quad b = B y,
\end{equation}
\eqref{eq:necessary} transforms into
\begin{equation}
  \label{eq:normal_form}
\big( A^\top \! A + \gamma \Theta(0) \big) w = A^\top b.
\end{equation}
By resorting to the Woodbury matrix identity, the solution to \eqref{eq:normal_form} can be expressed in a computationally efficient form
\begin{equation}
  \label{eq:LD_mean}
w = \Theta(0)^{-1} A^\top \big(\gamma \mathrm{I} + A \Theta(0)^{-1} A^\top \big)^{-1} b.
\end{equation}
Note that \eqref{eq:normal_form} and \eqref{eq:LD_mean} correspond to Tikhonov regularization with  $H^1(\Omega)$ penalization; the main reason for introducing this first algorithm as solving an approximate version of the necessary condition \eqref{eq:necessary} is allowing a direct comparison between the two algorithms.

The difference between the forward models \eqref{eq:lin_mod} and \eqref{eq:proj_lin_mod} manifests itself in the matrix $A$ via the symmetric factor $B$ in \eqref{eq:Amatrix}. For \eqref{eq:lin_mod}, one can simply set $B = C$, where $C$ is a Cholesky factor for $\Gamma_E^{-1}$, i.e., $\Gamma_E^{-1} = C^\top C$.  On the other hand, $B = C P$ is an appropriate choice for \eqref{eq:proj_lin_mod}.

\subsection{Total variation with lagged diffusivity iteration}
\label{sec:L-D}
The second algorithm employs the lagged diffusivity iteration \cite{Vogel96} to minimize \eqref{eq:Tikhonovk}, i.e., to solve the linearized EIT problem \eqref{eq:necessary} with the nonlinear TV prior \eqref{eq:TV}. The algorithm starts from an initial guess for the conductivity perturbation $w^{(0)} =0$, and like the first algorithm, computes a reconstruction, say, $w^{(1)}$ using \eqref{eq:LD_mean}. This reconstruction is then taken as a new initial value based on which one builds a new regularization matrix $\Theta(w^{(1)})$ and proceeds to find an updated reconstruction. This process is continued iteratively, and it is described by the equation
\begin{equation}
  \label{eq:LD_mean2}
w^{(j+1)} =  \Theta(w^{(j)})^{-1} A^\top \big(\gamma \mathrm{I} + A \Theta(w^{(j)})^{-1} A^\top \big)^{-1} b.
\end{equation}
The iteration \eqref{eq:LD_mean2} can be interpreted as computing successive MAP estimates with a sequence of Gaussian priors that approximate the TV prior \cite{Helin22}. By utilizing the convexity of the minimization target \eqref{eq:Tikhonovk}, it could also be shown that the iteration converges towards the minimizer of \eqref{eq:Tikhonovk} in our finite-dimensional setting; see, e.g., \cite{chan1999convergence, dobson1997convergence} as well as \cite[Theorem~4.2]{Pohjavirta21} for a proof in a closely related setting. In our numerical experiments, the lagged diffusivity iteration is continued for long enough to observe clear convergence (or divergence); see, e.g., \cite{Arridge13,Harhanen15,Helin22} for more sophisticated stopping criteria.

\begin{remark}
\label{remark:multilin}
Both reconstruction algorithms can also be applied to the original nonlinear problem via sequential linearizations, that is, by introducing an outer iteration that corresponds to computing the Jacobian $J_\sigma$ at the latest approximate solution and using it to form updated versions of the linearized equations \eqref{eq:lin_mod} and \eqref{eq:proj_lin_mod}. In case of the lagged diffusivity iteration, one must pay some attention to make sure that the TV penalization acts on $\sigma - \sigma_0$ and not on the increment to the solution between successive outer iterations (cf.~\cite[Remark~3.1]{Hyvonen24}), but otherwise the implementation of such nonlinear solvers is straightforward; see,~e.g.,~\cite{Harhanen15,Hyvonen24} for algorithms based on combining sequential linearizations and the lagged diffusivity iteration in EIT. However, we limit ourselves here to considering the linearized equations \eqref{eq:lin_mod} and \eqref{eq:proj_lin_mod} since they suffice for demonstrating the advantages of projecting by $P_\zeta(\sigma_0, \zeta_0)$ and $P_{\zeta,\phi}(\sigma_0, \zeta_0)$.
\end{remark}

\section{Experiments}
\label{sec:experiments}

This section experiments on using the projections $P_{\zeta}(\sigma_0, \zeta_0)$ and $P_{\zeta, \phi}(\sigma_0, \zeta_0)$ to form reconstructions from real-world data with a high level of uncertainty about the quality of the contacts at the electrodes. We use the same data sets as in Section~\ref{sec:Paction}, collected using the experimental setup described in Section~\ref{sec:experimental}. In particular, we consider the the {\em resistor test case} and {\em tape test case} that, respectively, correspond to using adjustable resistors in the electrode cables and covering parts of the electrodes with tape; see Table~\ref{tab:resistors_and_tapes} and Figure~\ref{fig:taped_electrodes}, where the left image shows the cylindrical inclusion used in the resistor test case and the right image depicts that for the tape test case. Consult Sections~\ref{sec:experimental} and \ref{sec:Paction} for more details.

\begin{figure}[t!]
\center{
  {\includegraphics[width=10cm]{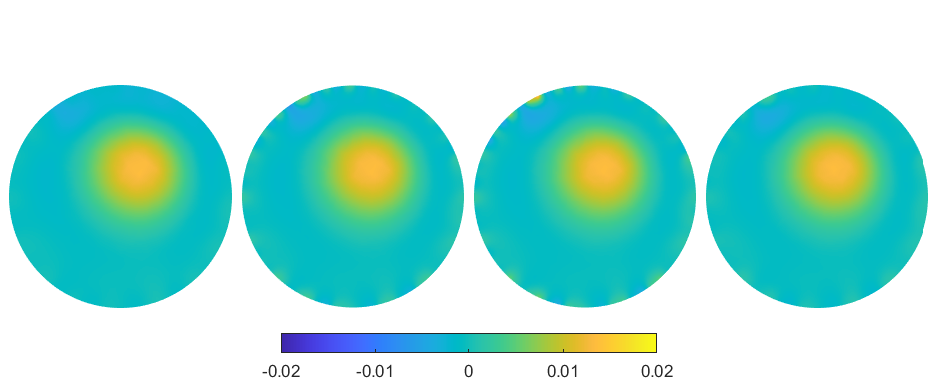}}
  {\includegraphics[width=10cm]{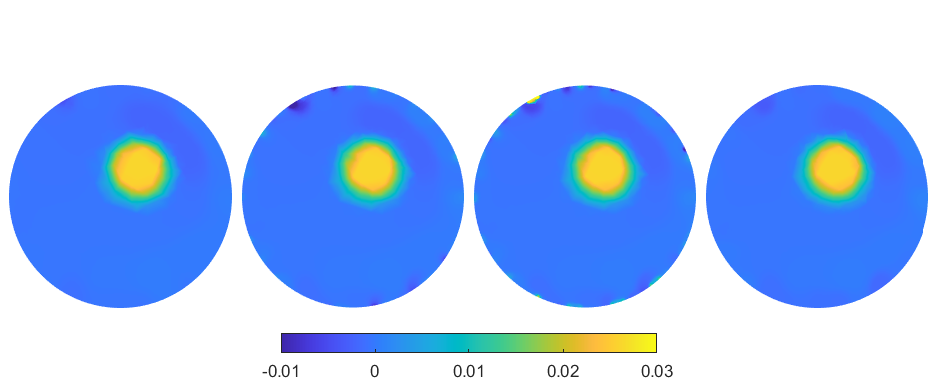}}
  }
\caption{Reference reconstructions for the left-hand setup in Figure~\ref{fig:taped_electrodes}, corresponding to the resistor test case,  without modifications to the electrode contacts. For both reconstructions horizontal slices at heights of 1\,cm, 2\,cm, 2.5\,cm and 3.5\,cm are displayed. Top: one-step algorithm. Bottom: lagged diffusivity algorithm.}
\label{fig:rec_noR}
\end{figure}

\begin{figure}[t]
\center{
  {\includegraphics[width=10cm]{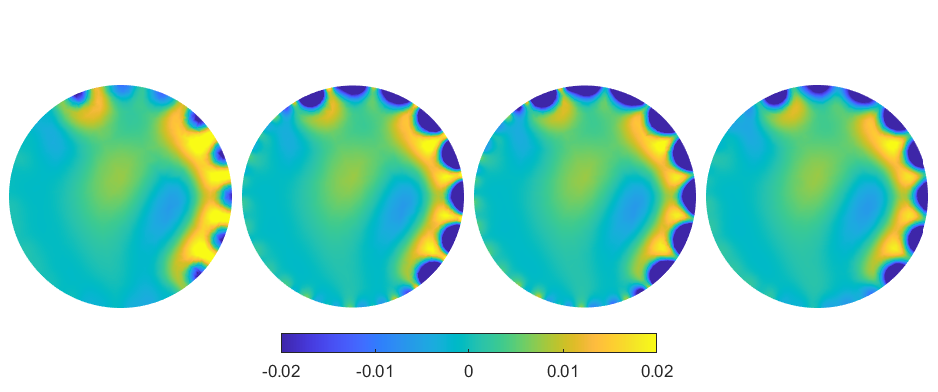}}
  {\includegraphics[width=10cm]{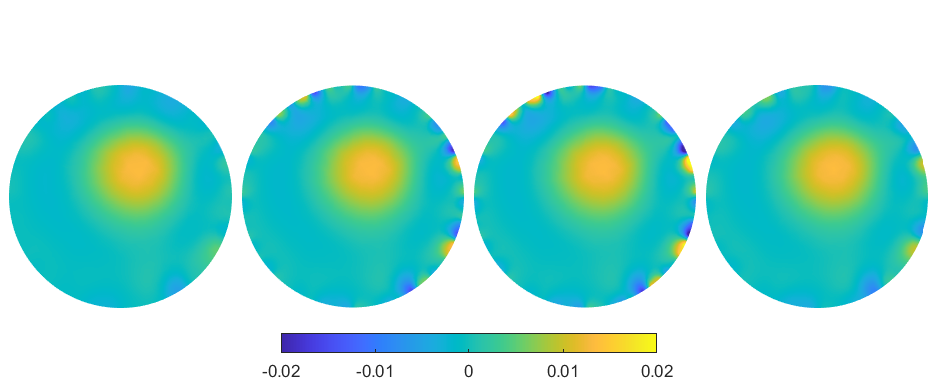}}
  {\includegraphics[width=10cm]{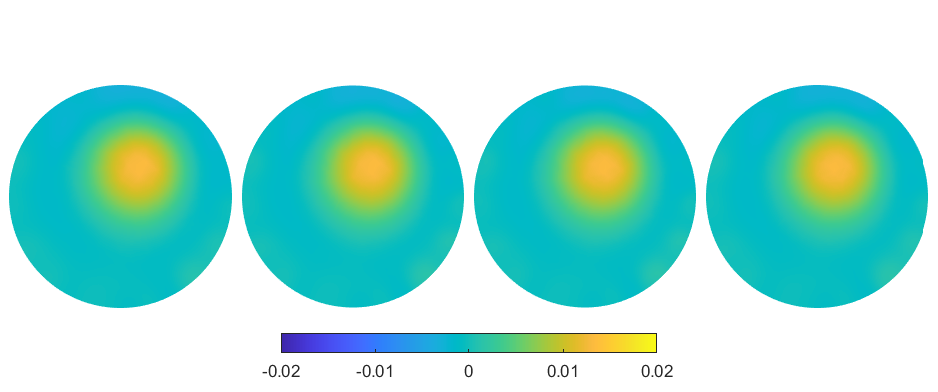}}
  }
\caption{Reconstructions computed by the one-step algorithm for the resistor test case considered on the left in Figure~\ref{fig:taped_electrodes}. For each reconstruction horizontal slices at heights of 1\,cm, 2\,cm, 2.5\,cm and 3.5\,cm are displayed. Top: No projections used. Middle: $P = P_{\zeta}(\sigma_0,\zeta_0)$ used for forming $B$ in \eqref{eq:Amatrix}. Bottom: $P = P_{\zeta,\phi}(\sigma_0,\zeta_0)$ used for forming $B$ in \eqref{eq:Amatrix}.}
\label{fig:rec_R_1step}
\end{figure}

\begin{figure}[t]
\center{
  {\includegraphics[width=10cm]{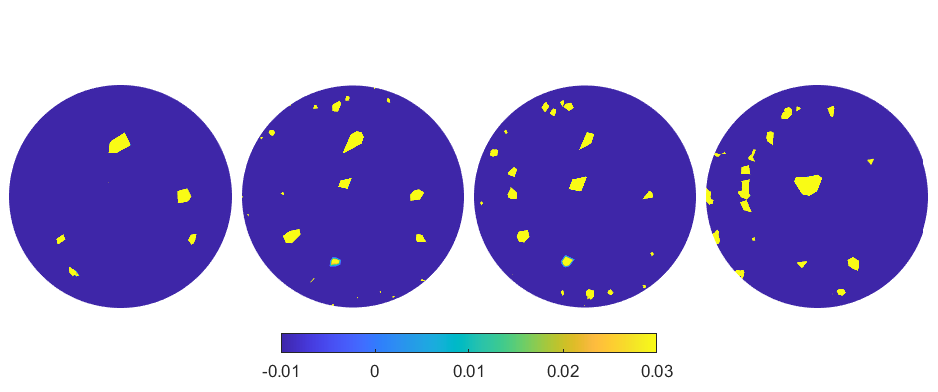}}
  {\includegraphics[width=10cm]{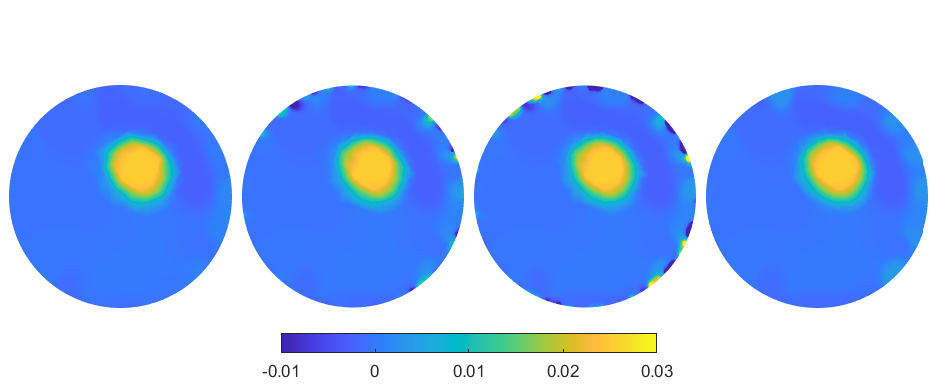}}
  {\includegraphics[width=10cm]{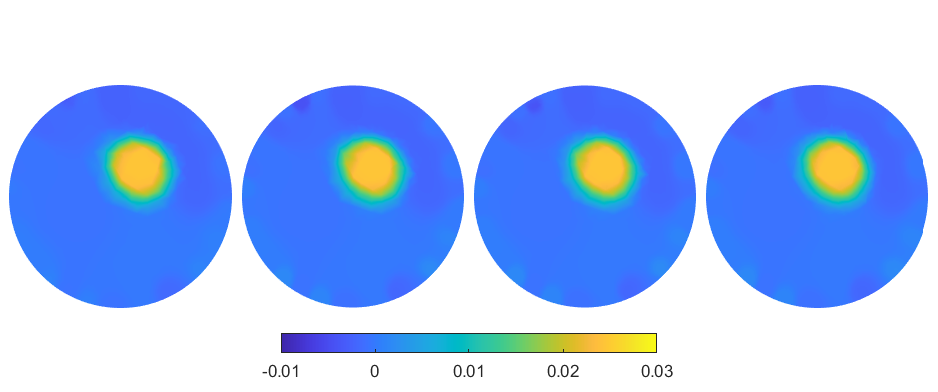}}
  }
\caption{Reconstructions computed by the lagged diffusivity algorithm for the resistor test case considered on the left in Figure~\ref{fig:taped_electrodes}. For each reconstruction horizontal slices at heights of 1\,cm, 2\,cm, 2.5\,cm and 3.5\,cm are displayed. Top: No projections used. Middle: $P = P_{\zeta}(\sigma_0,\zeta_0)$ used for forming $B$ in \eqref{eq:Amatrix}. Bottom: $P = P_{\zeta,\phi}(\sigma_0,\zeta_0)$ used for forming $B$ in \eqref{eq:Amatrix}.}
\label{fig:rec_R_10step}
\end{figure}

\begin{figure}[t]
\center{
  {\includegraphics[width=10cm]{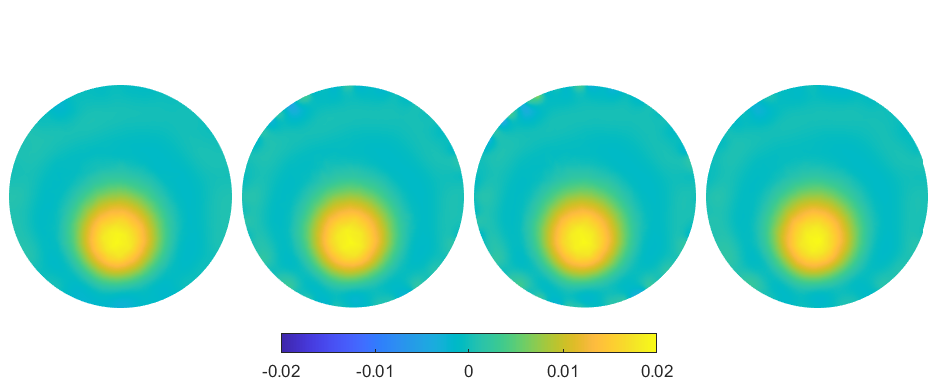}}
  {\includegraphics[width=10cm]{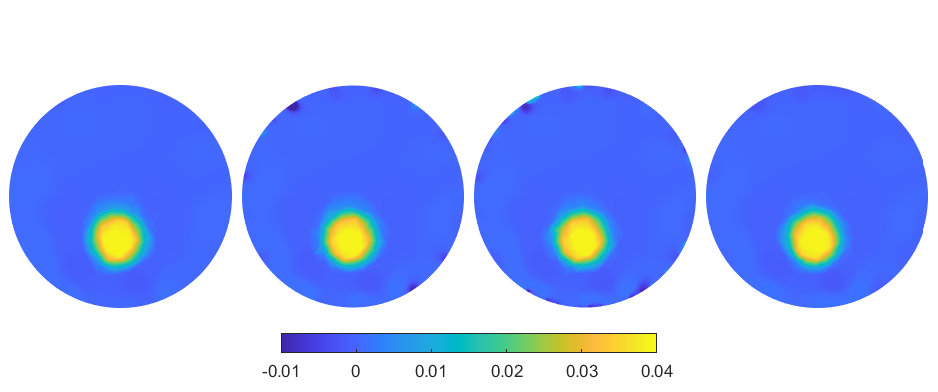}}
  }
\caption{Reference reconstructions for the right-hand setup in Figure~\ref{fig:taped_electrodes}, corresponding to the tape test case,  without modifications to the electrode contacts. For both reconstructions horizontal slices at heights of 1\,cm, 2\,cm, 2.5\,cm and 3.5\,cm are displayed. Top: one-step algorithm. Bottom: lagged diffusivity algorithm.}
\label{fig:rec_noTape}
\end{figure}

\begin{figure}[t]
\center{
  {\includegraphics[width=10cm]{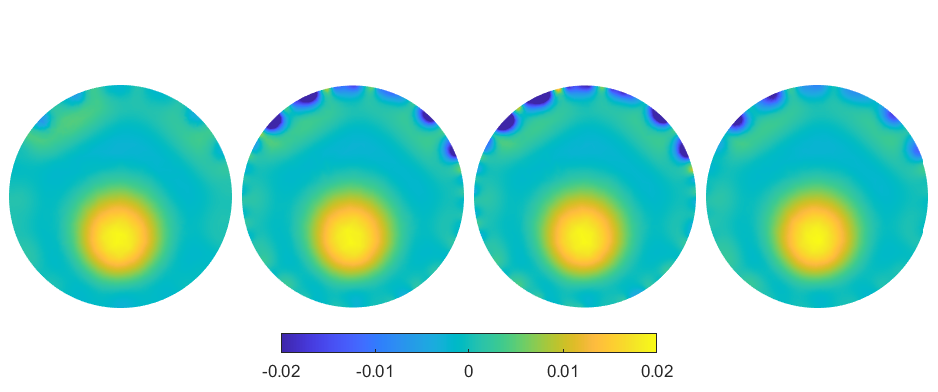}}
  {\includegraphics[width=10cm]{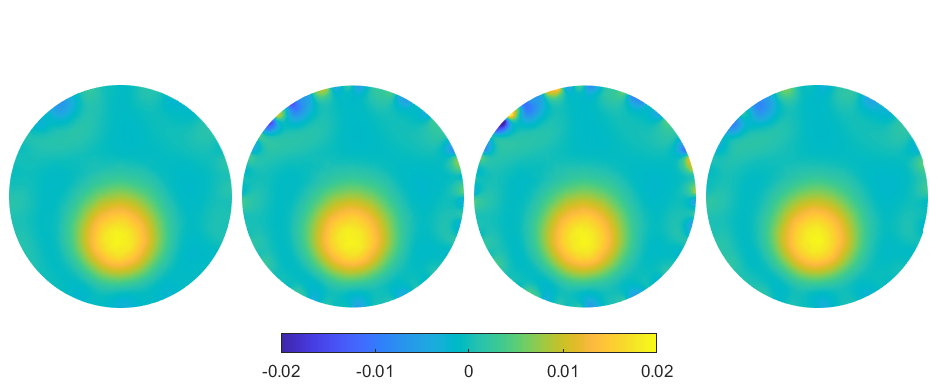}}
  {\includegraphics[width=10cm]{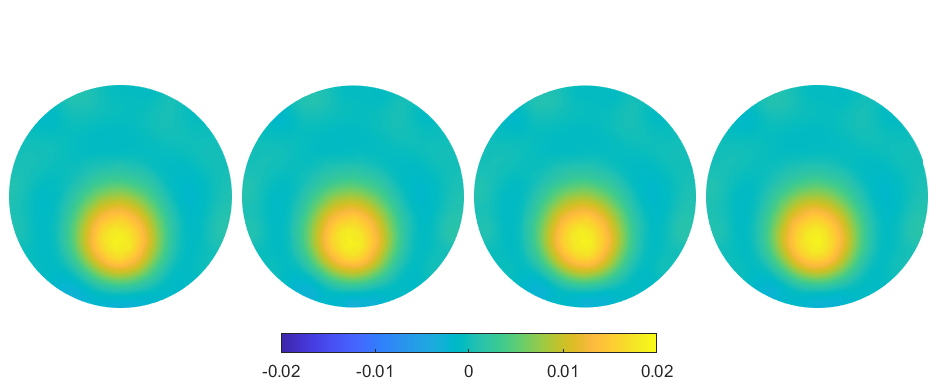}}
  }
\caption{Reconstructions computed by the one-step algorithm for the tape test case considered on the right in Figure~\ref{fig:taped_electrodes}. For each reconstruction horizontal slices at heights of 1\,cm, 2\,cm, 2.5\,cm and 3.5\,cm are displayed. Top: No projections used. Middle: $P = P_{\zeta}(\sigma_0,\zeta_0)$ used for forming $B$ in \eqref{eq:Amatrix}. Bottom: $P = P_{\zeta,\phi}(\sigma_0,\zeta_0)$ used for forming $B$ in \eqref{eq:Amatrix}.}
\label{fig:rec_maxTape_1step}
\end{figure}

\begin{figure}[t]
\center{
  {\includegraphics[width=10cm]{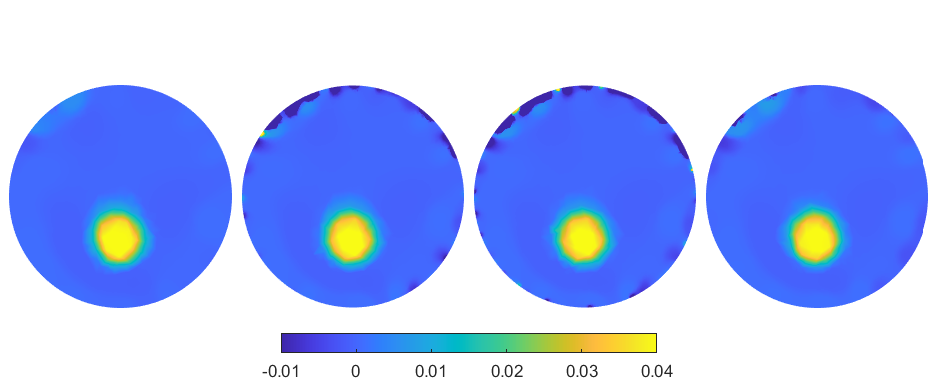}}
  {\includegraphics[width=10cm]{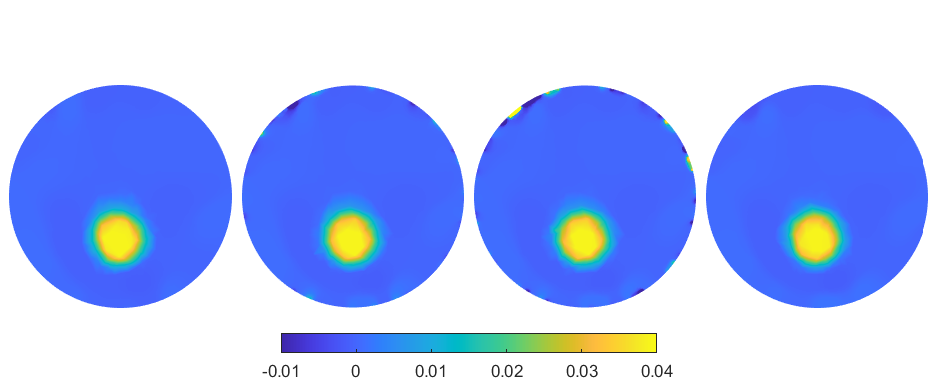}}
  {\includegraphics[width=10cm]{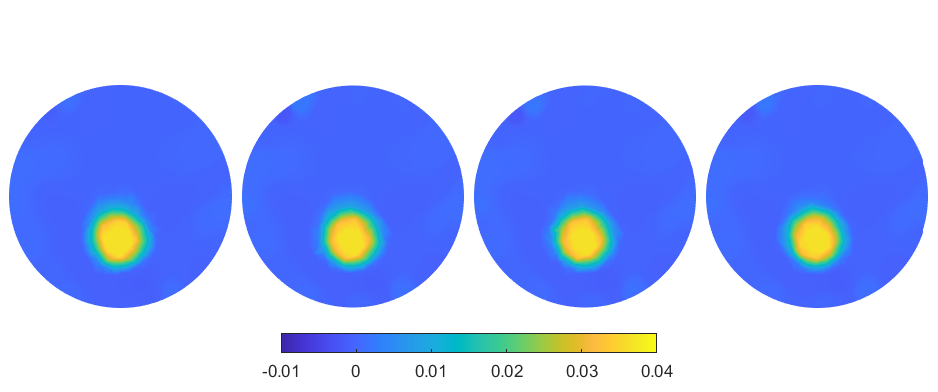}}
  }
\caption{Reconstructions computed by the lagged diffusivity algorithm for the tape test case considered on the right in Figure~\ref{fig:taped_electrodes}. For each reconstruction horizontal slices at heights of 1\,cm, 2\,cm, 2.5\,cm and 3.5\,cm are displayed. Top: No projections used. Middle: $P = P_{\zeta}(\sigma_0,\zeta_0)$ used for forming $B$ in \eqref{eq:Amatrix}. Bottom: $P = P_{\zeta,\phi}(\sigma_0,\zeta_0)$ used for forming $B$ in \eqref{eq:Amatrix}.}
\label{fig:rec_maxTape_10step}
\end{figure}

Both algorithms described in Section~\ref{sec:algorithms} are used to compute reconstructions; the results are presented in Figures~\ref{fig:rec_noR}--\ref{fig:rec_maxTape_10step}. For the one-step algorithm of Section~\ref{sec:one-step}, we choose the parameter value $\gamma = 10^{-2}$ in \eqref{eq:LD_mean}, and for the lagged diffusivity algorithm, we use $\gamma = 10^2$ in \eqref{eq:LD_mean2}. The values of the small free parameters $T$ and $\varepsilon$ are chosen, respectively, to be $10^{-6}$ and  the second smallest eigenvalue of $\Theta(w)$, but their precise values do not have any major effect on the conclusions of the numerical experiments. The covariance matrix for the measurement noise is defined by assuming that the noise components are mutually independent and have a common standard deviation that is $0.5\%$ of the maximum difference between potential values measured for the empty tank without artificially worsened contacts, resulting in $\Gamma_E = 1.6 \cdot 10^{-4}\, {\rm I}$ for the resistor test case and $\Gamma_E = 2.6 \cdot 10^{-4} \, {\rm I}$ for the tape test case. The lagged diffusivity iteration \eqref{eq:LD_mean2} is run for 10 iterations in all experiments, which is enough for sufficient convergence according to visual inspection. As in Section~\ref{sec:Paction}, the used background conductivity $\sigma_0$ corresponds to an experimentally measured value for the conductivity of the salt water in the tank, which is 0.0491\,S/m for the resistor test case and 0.0322\,S/m for the tape test case. The background peak contact conductivity values $\zeta_0$ are estimated based on measurements with the empty tank (cf.,~e.g.,~\cite[Section~4.3]{Hyvonen17b}), and they correspond to no extra resistors or tapes affecting the contacts. Recall that $\sigma_0$ and $\zeta_0$ are used for computing the reference measurements $\mathcal{U}(\sigma_0, \zeta_0)$ as well as the projections $P_{\zeta}(\sigma_0, \zeta_0)$ and $P_{\zeta, \phi}(\sigma_0, \zeta_0)$ for the algorithms of Section~\ref{sec:algorithms}, which means that the employed prior information on the contacts is highly inaccurate for both the resistor and the tape test cases. Recall also that the utilized projection enters the reconstruction algorithms via the symmectric matrix factor $B$ in \eqref{eq:Amatrix}. The values in the reconstruction images present the difference of the reconstructed conductivity from the constant background value $\sigma_0$.

Figure \ref{fig:rec_noR} illustrates the reference reconstructions for the left-hand setup in Figure~\ref{fig:taped_electrodes} when no resistors have been added to the cables, i.e., when the prior information on the peak values of the contact conductivity is relatively accurate. These reconstructions serve as a baseline for how accurate of a reconstruction the algorithms can provide when no projections are used and the contact conductivities are approximately known. In this case, both algorithms perform reasonably well in locating the highly conductive inclusion. As expected, due to the proper use of the TV prior, the reconstruction given by the lagged diffusivity algorithm looks slightly more accurate. However, both algorithms significantly underestimate the conductivity of the inclusion that was measured to be 4.73\,S/m, which is to be expected as the reconstructions are based on a single linearization and is also quite typical for EIT algorithms in general. Additionally, the reconstructions contain some small artifacts near the electrodes. These artifacts are presumably caused by inaccuracies in the contact conductivities: even though the peak values for the contact conductivity have been estimated based on empty tank measurements, these estimates can be somewhat inaccurate, we do not know the contact conductivity profiles on the electrodes, and the contacts may have slightly changed between the measurements for the empty tank and with an embedded cylinder. Despite the inaccuracies, these reconstructions provide good baseline information on the conductivity change in the water tank, and by using the projections $P_{\zeta}(\sigma_0, \zeta_0)$ and $P_{\zeta, \phi}(\sigma_0, \zeta_0)$ in the reconstruction algorithms, we hope to be able to achieve results of similar quality even when the contact conductivities are not known (and their estimation is not included in the algorithms).

Figure \ref{fig:rec_R_1step} demonstrates the effects of the projections on reconstructions computed using the one-step algorithm in the resistor test case. As can be seen from the reconstruction on the top row of Figure \ref{fig:rec_R_1step}, the resistors in the electrode cables add large errors to the reconstruction close to the manipulated electrodes when no projection is used. The algorithm interprets the resistors as more resistive regions near the electrodes, which almost completely overshadow the inclusion that we want to detect. By projecting with $P_\zeta(\sigma_0, \zeta_0)$, i.e., with respect to the contact conductivities, the conductivity fluctuations are mostly, but not completely removed from the reconstruction, and the location of the true inclusion is recovered; see the reconstruction on the middle row. As illustrated by the reconstruction on the bottom row of Figure~\ref{fig:rec_R_1step}, projecting also with respect to the electrode locations, i.e., employing $P = P_{\zeta,\phi}(\sigma_0, \zeta_0)$ when forming $B$ for \eqref{eq:Amatrix}, removes all remaining visible artifacts. In fact, this reconstruction is of better quality than the corresponding reference one in Figure~\ref{fig:rec_noR}.

Figure~\ref{fig:rec_R_10step} visualizes the results of an analogous experiment for the resistor test case as Figure~\ref{fig:rec_R_1step}, but with the reconstructions computed using the lagged diffusivity algorithm of Section~\ref{sec:L-D}. In this case, when no projections are used, the algorithm breaks down and does not produce a useful reconstruction. This happens because the algorithm starts by creating resistive regions near the electrodes with resistors, and the successive lagged diffusivity iterations then push these regions to become smaller in size with larger and larger conductivity fluctuations. Via fine-tuning the value of~$\gamma$, something more similar in quality to the top reconstruction in Figure \ref{fig:rec_R_1step} could be obtained, but we choose to keep the value of $\gamma$ constant between reconstructions in order allow a direct comparison. As with the one-step algorithm, projecting with $P_\zeta(\sigma_0, \zeta_0)$ removes most of these issues, only leaving small artifacts close to the affected electrodes. Including the projection with respect to electrode locations, i.e., using $P_{\zeta, \phi}(\sigma_0, \zeta_0)$, we obtain a reconstruction of good quality --- even better than the corresponding one in Figure~\ref{fig:rec_noR} --- that accurately locates the inhomogeneity and carries only minor artifacts in its background.

The corresponding results for the tape test case are documented in Figures~\ref{fig:rec_noTape}, \ref{fig:rec_maxTape_1step} and  \ref{fig:rec_maxTape_10step}, which are organized in the same way as Figures~\ref{fig:rec_noR}, \ref{fig:rec_R_1step} and \ref{fig:rec_R_10step}, respectively, for the resistor test case. This time, the reconstructions without employing either of the projections on the top rows of Figures~\ref{fig:rec_R_1step} and \ref{fig:rec_R_10step} are able to reveal the location of the inclusion, but the artifacts caused by the taped electrodes are in any case considerable. The reason for the better performance is probably three-fold: the considered inclusion is larger than the one in the resistor test case, it lies further away from the electrodes whose contacts have been worsened, and the tapes affect the overall contacts less than the resistors. Otherwise, the conclusions remain the same as in the resistor test case: utilizing the contact conductivity projection $P_\zeta(\sigma_0, \zeta_0)$ in the reconstruction algorithms significantly reduces the artifacts caused by mismodeling of the contact strengths, and resorting to the projection accounting also for the azimuthal angles of the electrodes, $P_\zeta(\sigma_0, \zeta_0)$, leads to reconstructions that are even better than the reference ones in Figure~\ref{fig:rec_noTape}.

\section{Concluding remarks}
\label{sec:conclusion}
This work demonstrated that projecting EIT measurements onto the orthogonal complement of the range of the Jacobian with respect to the contact conductivities and the electrode positions has potential to (i) reveal if an observed change in the measurements originates from changes in the electrode contacts or from a change in the internal conductivity of the imaged body and to (ii) allow reconstructing the internal conductivity without a need to simultaneously estimate the contacts. The former has a potential application in, e.g., monitoring of stroke patients in intensive care \cite{Toivanen21,Toivanen24}, and the latter can simplify reconstruction algorithms by enabling to ignore the estimation of contacts during the reconstruction process. What is more, according to our numerical tests, the range of the Jacobian with respect to the contact conductivity seems to be almost independent of the values of the contact and internal conductivities at which it is evaluated (although the same does not hold for the Jacobian itself), which means the utilized projection matrix can be evaluated even at a significantly inaccurate initial guesses for the two conductivities without compromising the above conclusions. In addition, the ability to project away the change in the measurements due to the electrode contacts does not seem to depend on the used current patterns: although the presented tests with experimental data employed current patterns for each of which only two electrodes actively participated in driving the current, our numerical simulations demonstrated the phenomenon also when (almost) all electrodes have an active role in feeding each of the considered current patterns (cf.~Section~\ref{sec:head}).  

Future lines of research include testing the presented projection idea in more complicated and practical measurement setups, including online monitoring of stroke to deduce if the changes in the measurements are due to (secondary) bleeding in the brain or changes in the contacts. Another aim is to establish a better theoretical understanding on why and when one can utilize a simple linear projection to infer whether the changes in the measurements are caused by changes in internal or contact conductivities.

\bibliographystyle{acm}
\bibliography{contcondproj-refs.bib}
\end{document}